\documentclass[journal]{IEEEtran}  

\usepackage{graphics} 
\usepackage{epsfig} 
\usepackage{mathptmx} 
\usepackage{times} 
\usepackage{amsmath} 
\usepackage{amssymb}  
\usepackage{amsfonts}
\usepackage{diagbox}
\usepackage{amsthm}

\usepackage{algorithmic}
\usepackage{enumerate} 
\usepackage[usenames, dvipsnames]{color}
\usepackage{booktabs}
\usepackage{rotating}
\usepackage{setspace}
\usepackage{algorithm}
\usepackage[margins]{trackchanges}

\title{Robust model predictive control for large-scale distributed parameter systems under uncertainty}
\author{Min Tao, Ioannis Zacharopoulos and Constantinos Theodoropoulos
\thanks{M. Tao, I. Zacharopoulos and C. Theodoropoulos are with the Department of Chemical Engineering and Analytical Science, University of Manchester, Manchester M13 9PL, U.K. (e-mail: min.tao@sheffield.ac.uk;
ioannis.zacharopoulos@manchester.ac.uk;
k.theodoropoulos@manchester.ac.uk} 
}

\begin{document}
\maketitle
\begin{abstract}
Control of nonlinear distributed parameter systems (DPS) under uncertainty is a meaningful task for many industrial processes. However, both intrinsic uncertainty and high dimensionality of DPS require intensive computations, while non-convexity of nonlinear systems can inhibit the computation of global optima during the control procedure.  In this work, polynomial chaos expansion (PCE) was used to account for the uncertainties in quantities of interest through a systematic data collection from the high-fidelity simulator. Then the proper orthogonal decomposition (POD) method was adopted to project the high-dimensional nonlinear dynamics of the computed statistical moments/bounds onto a low-dimensional subspace, %
where recurrent neural networks (RNNs) were subsequently built to capture the reduced dynamics. 
Finally, the reduced RNNs based model predictive control (MPC) 
would generate a set of sequential optimisation problems, of which near global optima could be computed through the mixed integer linear programming (MILP) reformulation techniques and advanced MILP solver.
The effectiveness of the proposed framework is demonstrated through two case studies: a chemical tubular reactor and a cell-immobilisation packed-bed bioreactor for the bioproduction of succinic acid.
\end{abstract}

\section{Introduction}
\label{S:1}
Spatial-temporal distributed parameter systems (DPS) exist widely in engineering practice  \cite{ahmed1987distributed}, e.g in chemical \cite{tao2016hybrid}, biochemical \cite{park2015integration} and mechanical engineering \cite{yang2008simulation}. 
Complex DPS usually exhibit uncertainty due to inherent stochastic and/or incomplete knowledge of processes \cite{sullivan2015introduction}, which lead to substantial model-plant mismatch. Efficient control strategies for large-scale complex DPS under uncertainty could speed up process production and ensure process safety \cite{stephanopoulos1996intelligent}.
Model predictive control (MPC), a popular advanced control method for multivariate plants with process constraints, reformulates the original optimal control problem (OCP) into a finite sequence of dynamic optimisation problems at each sampling time to obtain the corresponding control actions \cite{rawlings2017model}. Feedback is introduced to this procedure through updating the states of nonlinear dynamic models. In spite of both improved optimisation algorithms and fast hardware, the computational demand on MPC is still high for large-scale DPS problems. The high dimensionality of the discretised DPS results in slow computational speeds, which is impractical for real-time control, while global optimisation of the non-convex MPC sub-problems is usually computationally intractable. In addition, black-box characteristics of high-fidelity commercial simulators \cite{multiphysics1998introduction, fluent2015ansys} prevent the direct utilisation of model-based computational techniques. Furthermore, dynamic model predictions are often significantly affected by parametric uncertainties \cite{nagy2007distributional}, which may lead to wrong decision-making from MPC controllers. Therefore, efficient uncertainty quantification methods are of great importance for developing MPC strategies for large-scale DPS under parametric uncertainty. Uncertainty quantification procedures for DPS however, often require a large number of repeated computationally expensive evaluations. 
Thus, control of DPS under uncertainty is an important practical engineering challenge.

%

Uncertainty quantification (UQ) aims to measure the impact of uncertainty on quantities of interest \cite{sullivan2015introduction}. The direct Monte Carlo (MC) sampling method is typically utilised to complete UQ tasks for generalised complex problems, generating a large number of realisations to accurately approximate the uncertainty distributions. Compared with the expensive MC method, both efficient sampling methods \cite{dunn2011uncertainty} and lower order model-based methods \cite{nagy2007distributional} are more powerful in terms of computational costs for complex large-scale systems. Efficient sampling methods, such as Latin hypercube sampling \cite{dunn2011uncertainty} and sparse grid methods \cite{nobile2008sparse,nobile2008anisotropic,lin2010uncertainty}, only require a few representative samples to be propagated utilising the original systems, which can greatly reduce computational demand. However, less samples would lead to lower computational accuracy for quantitative tasks. Lower order models on the other hand, including popular power series expansions (PSE) \cite{nagy2003worst,nagy2004open} and polynomial chaos expansions (PCE) \cite{eldred2009recent}, can lighten the computational load by replacing the original expensive computational models \cite{nagy2007distributional,bradford2019output}, which has been commonly employed to estimate the statistics for large-scale complex systems \cite{chaffart2017robust,chaffart2018robust}. Computational case studies have shown that the PCE method is faster and more accurate than the PSE method for large-scale thin film formation and heterogeneous catalytic flow problems \cite{kimaev2017comparison,chaffart2017robust}. In this work, fast PCE-based uncertainty propagation method was employed to address parametric uncertainty for large-scale complex systems. Nevertheless, both high dimensionality and non-convexity issues are still computational barriers for applying MPC to large-scale distributed systems. 

Model order reduction techniques are the most efficient methodologies to address the high dimensional issues for spatial-temporal distributed systems \cite{schilders2008model}. Together with Galerkin method \cite{theodoropoulou1998model,xie2011off} or artificial neural network (ANN) surrogate models \cite{xie2015data}, the proper orthogonal decomposition (POD) approach, also named as Karhunen-Loeve decomposition, can generate accurate low-dimensional models that can be efficiently used to perform control and optimisation tasks for large-scale distributed systems. Additionally, inertial manifolds or approximate inertial manifolds have been employed to construct stable controllers for dissipative partial differential equation (PDE) systems \cite{christofides1997finite,christofides1998robust}. Moreover,  \textit{equation-free} methodologies offer another effective model reduction approach for optimisation and control purposes \cite{kevrekidis2003equation,luna2005input,bonis2012model,petsagkourakis2018reduced}. Equation-free methods utilise input/output data to compute dominant system eigendirections  that {\it drive} the system dynamics as well as low-dimensional gradients that can accelerate computational evaluations for large-scale dissipative systems. A detailed discussion about model reduction based optimisation and control methodologies can be found in \cite{theodoropoulos2011optimisation}. 

To overcome the non-convexity issues for large-scale nonlinear optimisation problems, both stochastic and deterministic methods can be used. Stochastic optimisation methods \cite{kirkpatrick1983optimization,chambers2019practical} can perform global searches across the design space to avoid multiple local optima, but they can not guarantee optimal conditions. Deterministic global optimisation methods typically implement branch and bound procedures or their extensions to narrow the gap between low and upper bounds \cite{floudas2005global}, and are often extremely expensive for large-scale problems but they can guarantee global optimality conditions.  

The aim of this work is to construct a robust model predictive control framework for large-scale distributed systems under uncertainty. Firstly, polynomial chaos expansion approach is used to compute statistics for the quantities of interest of DPS. Then, a combination of POD and recurrent neural networks (RNNs) is employed to capture the nonlinear dynamics of the calculated statistical moments and/or bounds. The resulting reduced ANN surrogate models can accurately represent the dynamics of the (black-box) large-scale systems and are efficiently implemented for MPC. In addition, the resulting optimisation sub-problems are globally optimised using advanced global optimisation solvers \cite{misener2014antigone,tawarmalani2005polyhedral,rehfeldt2018scip}.

Since the high complexity of ANN structures leads to intractable computational problems for advanced global solvers, most previous studies of optimising surrogate ANN-based models focus on local optimisation \cite{henao2011surrogate} and/or small problems \cite{smith2012cfd}. A reduced-space global optimisation algorithm \cite{schweidtmann2019deterministic} for ANN-based models was proposed to reduce computational costs and to allow an efficient nonlinear MPC formulation. \cite{doncevic2020deterministic}.  
The reduced-space global optimisation algorithm focuses more on the online reformulation techniques within global optimisation procedures for general ANN models. .
In our previous work \cite{tao2021model}, offline reformulation strategies were exploited to build suitable ANN structures and the related mathematical formulations to improve the computational performances. Specifically, principal component analysis and deep rectifier neural networks were employed to build accurate but relatively simple ANN models for distributed parameter systems, which could be then formulated into mixed integer linear programming (MILP) problems and globally solved by advanced MILP solver CPLEX \cite{ilog2006cplex}. Here, offline reformulations techniques, POD and ReLU based RNN models were employed to describe the original high-dimensional nonlinear dynamics of statistical moments and/or bounds for the quantities of interest. Then the dynamic optimisation of distributed parameter systems could be reformulated into MILP problems solved by CPLEX, which could significantly reduce computational costs and capture the approximate global optima.
Global NMPC studies have been performed in previous literature \cite{degachi2015global,katz2020integrating}, employing advanced global solvers and mixed integer programming strategies. The normalised multi-parametric disaggregation technique was utilised to compute upper and lower bounds in a spatial branch and bound formulation for global  NMPC with multiple operating conditions \cite{wang2017globally}. Also, rigorous nonlinear MPC was reformulated into a sequence of mixed integer nonlinear programming problems, which were then solved globally \cite{long2004globally}. However, all these works are limited to small and or medium size problems. To the best of our knowledge, global NMPC techniques have not been utilised for controlling large-scale DPS under uncertainty. In addition, no other work has been reported using a combination of double model reduction involving POD and RNN for high dimensional nonlinear systems, combined with PCE to address parametric uncertainty. The novelty of this work is to provide a PCE-POD-RNN based robust NMPC strategy for large-scale distributed systems under uncertainty, where the resulting surrogate model based optimisation sub-problems are globally solved. The performance of the proposed computational framework is validated via a receding horizon NMPC formulation for a chemical tubular reactor and a cell-immobilisation packed-bed biochemical reactor.

The rest of the paper is organized as follows. In Section 2, the robust NMPC strategy framework is proposed. Furthermore, the detailed theoretical basis and implementation are provided. In Section 3, the model based control framework is verified by practical chemical and biochemical reactors. In Section 4, conclusions and further applications are discussed. 

\section{PCE-POD-RNN based robust NMPC Methodology}
In this part, the PCE-POD-RNN based robust nonlinear MPC strategy would be introduced. Firstly, the general optimal control problem is formulated. Then detailed polynomial chaos expansion, proper orthogonal decomposition and recurrent neural network parts are discussed, respectively. Finally, the general robust nonlinear model predictive control methodology is illustrated.
\subsection{Problem formulation}
The general optimal control problem for PDE-based 
distributed parameter systems with parametric uncertainty:
  \begin{equation}\label{eq1}
  \begin{aligned}
\min_{\boldsymbol{\lambda}} \quad & G(\boldsymbol {y},\boldsymbol{P},\boldsymbol{\lambda})  \\  
s.t.  \quad  \frac{\partial  \boldsymbol{y}}{\partial \tau} &=
D\{\frac{\partial  \boldsymbol{y}}{\partial \boldsymbol{x}}, \frac{\partial^{2} \boldsymbol {y}}{\partial \boldsymbol{x}^{2}},...,
\frac{\partial^{n} \boldsymbol {y}} {\partial \boldsymbol{x}^{n}},\boldsymbol{P}\}+g_{stm}( \boldsymbol{y},\boldsymbol{\lambda}) \\
&\left. A \{\frac{\partial  \boldsymbol{y}}{\partial \boldsymbol{x}},\frac{\partial^{2} \boldsymbol {y}}{\partial \boldsymbol{x}^{2}},...,\frac{\partial^{n} \boldsymbol {y}}{\partial \boldsymbol{x}^{n}},\boldsymbol{\lambda}\}\right|_{\boldsymbol{x}=\Omega}=h_{bds}( \boldsymbol{y},\boldsymbol{\lambda}) \\
&\left. \boldsymbol{y}\right|_{\tau=0} = \boldsymbol{y}_{0}\\
& g_{cons}(\boldsymbol {y},\boldsymbol{P},\boldsymbol{\lambda}) \leq 0
\end{aligned}
\end{equation}
where $G(\boldsymbol {y},\boldsymbol{P},\boldsymbol{\lambda}) : \mathbb{R}^{N_y} \times \mathbb{R}^{N_p}\times \mathbb{R}^{N_{\lambda}} \to \mathbb{R}$  denotes the objective cost function for the OCP, and the constraints are the PDE-based system equations with corresponding boundary and initial conditions. $\tau \in \mathbb{R}$  is the time dimension while $\boldsymbol{x}\in \mathbb{R}^{N_{\boldsymbol{x}}}$ are the space dimensions, $\boldsymbol{y}(t,\boldsymbol{x}) : \mathbb{R}\times \mathbb{R}^{N_{\boldsymbol{x}}} \to \mathbb{R}^{N_y}$ are the state variables, $\boldsymbol{P} \in \mathbb{R}^{N_p}$  the uncertain parameters, $\boldsymbol{\lambda}(\tau) : \mathbb{R} \to \mathbb{R}^{N_{\lambda}}$  are the manipulated variables, $D$ is the differential operator, $g_{stm}( \boldsymbol{y},\boldsymbol{\lambda}): \mathbb{R}^{N_y} \times \mathbb{R}^{N_{\lambda}} \to \mathbb{R}^{N_y} $ are the nonlinear parts of PDEs, $A$ is the operator on the boundary conditions, $\Omega$ is the boundary, $h_{bds}( \boldsymbol{y},\boldsymbol{\lambda}) : \mathbb{R}^{N_y} \times \mathbb{R}^{N_{\lambda}} \to \mathbb{R}^{N_y}$ are the function values for the boundary conditions, $\boldsymbol{y}_{0}\in \mathbb{R}^{N_y}$ are the initial values and $g_{cons}(\boldsymbol {y},\boldsymbol{P},\boldsymbol{\lambda}): \mathbb{R}^{N_y} \times \mathbb{R}^{N_p}\times \mathbb{R}^{N_{\lambda}} \to \mathbb{R}^{N_{con}}$ are the $N_{con}$ general constraints for state variables, manipulated variables and uncertain parameters. Here, we assumed that the uncertain parameters $\boldsymbol{P} \in \mathbb{R}^{N_p}$ were time-invariant.  

In general, the analytical solutions for the above Eq.(\ref{eq1}) are unavailable. Therefore, computing the numerical solutions is more practical. If we transform the above continuous dynamic systems into discrete ones, the discretised dynamic systems can be formulated as below:
\begin{equation}\label{eq2}
\begin{aligned}
\ \boldsymbol{y'}_{t+1} &=f_{dis}(\boldsymbol{y'}_{t},\boldsymbol{\lambda'}_{t},\boldsymbol{P})+ \boldsymbol{w}_{t}\\
\ \boldsymbol{z}_{t+1} &=h_{dis}(\boldsymbol{y'}_{t+1},\boldsymbol{\lambda'}_{t},\boldsymbol{P})+ \boldsymbol{v}_{t}\\
\ \boldsymbol{y'}_{t=0}&=\boldsymbol{y}_{0} \\
\end{aligned}
\end{equation}where $\boldsymbol{y'}_{t}\in \mathbb{R}^{N_{y'}}$ the state variables at discretised time step $t \in\mathbb{N}$, 
$\boldsymbol{\lambda'}_{t}\in \mathbb{R}^{N_{\lambda'}}$ the manipulated values at discretised time step $t$, $\boldsymbol{z}_{t}\in \mathbb{R}^{N_{z}} $ the measurement values at discretised time step $t$.  $ \boldsymbol{y'}_{0}$ are the discretised initial values of $\boldsymbol{y}_{0} $. $f_{dis}: \mathbb{R}^{N_{y'} }\times \mathbb{R}^{N_{\lambda'}}\times \mathbb{R}^{N_{p}} \to \mathbb{R}^{N_{y'}}$ and $h_{dis}: \mathbb{R}^{N_{y'} }\times \mathbb{R}^{N_{\lambda'}}\times \mathbb{R}^{N_{p}} \to \mathbb{R}^{N_{z}}$ denote the discrete time nonlinear dynamic systems and output equations, respectively while  $\boldsymbol{w}_{t}\in \mathbb{R}^{N_{w}}$ and $\boldsymbol{v}_{t}\in \mathbb{R}^{N_{v}}$ represent the time-varying additional noises to state variables and measurements, respectively, which follow zero-mean normal distributions. 

Then the time-space infinite dimensional original problem Eq.(\ref{eq1}) can be reduced as a time-space finite dimensional large-scale stochastic (due to parametric uncertainty $\boldsymbol{P}$) 
nonlinear programming (NLP) problem:
  \begin{equation}\label{eq3}
  \begin{aligned}
\min_{\boldsymbol{\lambda'}_{t}} \quad & G'(\boldsymbol {y'}_{t},\boldsymbol{z}_{t},\boldsymbol{P},\boldsymbol{\lambda'}_{t})  \\   
s.t.  \ \boldsymbol{y'}_{t+1} &=f_{dis}(\boldsymbol{y'}_{t},\boldsymbol{\lambda'}_{t},\boldsymbol{P})+ \boldsymbol{w}_{t}, t=0,1,2....k\\
\ \boldsymbol{z}_{t+1} &=h_{dis}(\boldsymbol{y'}_{t+1},\boldsymbol{\lambda'}_{t},,\boldsymbol{P})+ \boldsymbol{v}_{t},t=0,1,2....k\\
\ \boldsymbol{y'}_{0}&=\boldsymbol{y}_{0} \\
 g'_{cons} &(\boldsymbol {y'_{t}},\boldsymbol{P},\boldsymbol{z}_{t},\boldsymbol{\lambda'}_{t}) \leq 0, t=0,1,2....k\\
\end{aligned}
\end{equation}
where $G'$ is the discretised objective cost function, $g'_{cons}$ denotes the general constraints for the discrete time state variables, manipulated variables, output measurements and uncertain parameters, and $k$ is the number of time horizons. Both the cost function $G'$ and the constraints $g'_{cons}$ have a stochastic representation due to the general parametric uncertainty $\boldsymbol{P}$. 
In the following parts, a combination of PCE and POD-RNN techniques is employed to deal with the large-scale stochastic nonlinear programming Eq.(\ref{eq3}).   

\subsection{Polynomial chaos expansion}
To address parametric uncertainty, probabilistic approaches describe it by employing a probability density function (PDF) \cite{stanton2000probability}. PCE, one of the most efficient UQ methods, uses only a few system samples to construct accurate stochastic surrogate models. The key idea of PCE is to represent an arbitrary random variable $g$ with finite second-order moments as a function of random variables $\boldsymbol{\theta}$ \cite{xiu2002wiener}:
\begin{equation}\label{eq4}
\begin{aligned}
g(\boldsymbol{\theta}) &=\sum a_{b}{\Theta}_{b}(\boldsymbol{\theta})\\
{\Theta}_{b}(\boldsymbol{\theta}) &=\Pi_{i=1}^{{N_{\theta}}}{\Theta}_{b_{i}}(\theta_{i})\\
\end{aligned}
\end{equation}
where $\boldsymbol{\theta} \in \mathbb{R}^{N_{\theta}}$ are random variables such as the random parameters $\boldsymbol{P}$ in this work, $g$ is quantities of interest such as the state variables and output values in this work; ${\Theta}_b : \mathbb{R}^{N_{\theta}} \to \mathbb{R}$ are multivariate orthogonal polynomials from tensor products of univariate polynomials $ {\Theta}_{b_{i}}:\mathbb{R}\to \mathbb{R}$  and $a_{b} \in \mathbb{R}$ are the corresponding coefficients, $b \in \mathbb{N}^{N_{\theta}}$ are multidimensional summation indices and  $b_{i}$ denotes the degree of each univariate polynomial ${\Theta}_{b_{i}}(\theta_{i})$ of $\theta_{i}$.

For generalised polynomial chaos, the choice of orthogonal polynomials significantly depends on the types of the probabilistic distributions of random variables $\boldsymbol{\theta}$. For example, Hermite polynomials would commonly be chosen for normal distributions. The equation below gives the expression for univariate Hermite polynomials ${\Theta}^{H}_{b_{i}}$ with respect to the standard Gaussian distribution $\theta^{g}_{i}$.
\begin{equation}\label{eq5}
{\Theta}^{H}_{b_{i}}(\theta^{g}_{i})=(-1)^{b_{i}}exp(\frac{1}{2}{\theta^{g}}^{2}_{i})\frac{d^{b_{i}}}{d{{\theta^{g}}^{b_{i}}_{i}}}exp(-\frac{1}{2}{\theta^{g}}^{2}_{i})\\
\end{equation}
One of the most important properties of multivariate orthogonal polynomials is orthogonality, i.e. generalised polynomial chaos terms are orthogonal to each other:
\begin{equation}\label{eq6}
<{\Theta}_{b1}(\boldsymbol{\theta}),{\Theta}_{b2}(\boldsymbol{\theta})>=  
\begin{cases}
<{\Theta}^{2}_{b1}(\boldsymbol{\theta})> & ,  b1 =b2\\
      0 & , b1 \neq b2\\
 \end{cases} 
\end{equation}
where $<{\Theta}^{2}_{b1}(\boldsymbol{\theta})>$ are often known as constants, whose values depend on the chosen polynomial family and multidimensional summation indice$b1$.

If the inner product $<,>$ of any two random variables $g_{2}(\boldsymbol{\theta}),g_{1}(\boldsymbol{\theta})$
are defined in corresponding probability space:
\begin{equation}\label{eq7}
<g_{1}(\boldsymbol{\theta}),g_{2}(\boldsymbol{\theta})>=\int_{\omega}g_{1}(\boldsymbol{\theta})g_{2}(\boldsymbol{\theta})\pi(\boldsymbol{\theta})d\boldsymbol{\theta}\\
\end{equation}where $\omega$ is the integral space of inner product and $\pi(\boldsymbol{\theta})$ is the probability density function, then the inner product between arbitrary random variable $g(\boldsymbol{\theta})$ and any multivariate orthogonal polynomial chaos ${\Theta}_{b1}(\boldsymbol{\theta})$ can be computed:
\begin{equation}\label{eq8}
\begin{aligned}
<{\Theta}_{b1}(\boldsymbol{\theta}),g(\boldsymbol{\theta})>&=\int_{\omega}{\Theta}_{b1}(\boldsymbol{\theta})g(\boldsymbol{\theta})\pi(\boldsymbol{\theta})d\boldsymbol{\theta}\\
&=\int_{\omega}{\Theta}_{b1}(\boldsymbol{\theta})\sum a_{b}{\Theta}_{b}(\boldsymbol{\theta})\pi(\boldsymbol{\theta})d\boldsymbol{\theta}\\
&=\sum a_{b}\int_{\omega}{\Theta}_{b1}(\boldsymbol{\theta}){\Theta}_{b}(\boldsymbol{\theta})\pi(\boldsymbol{\theta})d\boldsymbol{\theta}\\
&=\sum a_{b}<{\Theta}_{b1}(\boldsymbol{\theta}),{\Theta}_{b}(\boldsymbol{\theta})>\\
&= a_{b1}{<{\Theta}^{2}_{b1}(\boldsymbol{\theta})> }
\end{aligned}
\end{equation}
According to Eqs.(\ref{eq7}-\ref{eq8}), any coefficient of the generalised polynomial chaos terms can be computed as follows:
\begin{equation}\label{eq9}
\begin{aligned}
 a_{b1}&=\frac{<{\Theta}_{b1}(\boldsymbol{\theta}),g(\boldsymbol{\theta})>}{<{\Theta}^{2}_{b1}(\boldsymbol{\theta})> }\\
 &=\frac{\int_{\omega}{\Theta}_{b1}(\boldsymbol{\theta})g(\boldsymbol{\theta})\pi(\boldsymbol{\theta})d\boldsymbol{\theta}}{<{\Theta}^{2}_{b1}(\boldsymbol{\theta})> }
 \end{aligned}
\end{equation}
Here we adopt non-intrusive projection and quadrature methods to calculate the integral term as Eq.({\ref{eq9_1}}) due to the black-box characteristics of nonlinear dynamic systems.
\begin{equation}\label{eq9_1}
\begin{aligned}
\int_{\omega}{\Theta}_{b1}(\boldsymbol{\theta})g(\boldsymbol{\theta})\pi(\boldsymbol{\theta})d\boldsymbol{\theta}\approx\sum_{kk=1}^{N_{kk}} w_{kk}{\Theta}_{b1}(\boldsymbol{\theta_{kk}})g(\boldsymbol{\theta_{kk}})
 \end{aligned}
\end{equation}where $N_{kk}$ is the number of quadrature points,  $\theta_{kk}$ and $w_{kk}$ are the sampling points and the corresponding weights , respectively, according to the quadrature rules.

Then, the PCE-based stochastic surrogate models are tractable and can be utilised to estimate the probabilistic distributions of state variables and output values. The mean value $\mu_{\boldsymbol{g}}$ and covariance $\sigma_{\boldsymbol{g}}$ of quantity of interest $g(\boldsymbol{\theta})$
can be calculated as below \cite{paulson2014fast}:
\begin{equation}\label{eq10}
\begin{aligned}
g(\boldsymbol{\theta}) &\approx \sum_{b=\boldsymbol{0}}^{\boldsymbol{L}} a_{b}{\Theta}_{b}(\boldsymbol{\theta})\\
\mu_{\boldsymbol{g}}&=\mathbb{E}(g(\boldsymbol{\theta}))\approx a_{\boldsymbol{0}}\\
\sigma_{\boldsymbol{g}}&=\boldsymbol{Var}(g(\boldsymbol{\theta}))\approx\sum_{b=\boldsymbol{1}}^{\boldsymbol{L}} a_{b}^{2}<{\Theta}_{b}^2(\boldsymbol{\theta})>\\
\end{aligned}
\end{equation}where $\boldsymbol{L}$ denotes the tensor product of the truncated order of the arbitrary random variable $g(\boldsymbol{\theta})$. 

The specific probability limits, such as the worst bounds in low confidence levels under the defined probability $\mathbb{P}$, could be evaluated through the PDFs obtained for $g(\boldsymbol{\theta})$ as follows:
\begin{equation}\label{eq11}
\begin{aligned}
g^p(\boldsymbol{\theta}) &= F^{-1}_{cdf}( \mathbb{P})\\
\mathbb{P}(g(\boldsymbol{\theta})<g^p(\boldsymbol{\theta}) )& =F_{cdf}(g^p(\boldsymbol{\theta}) )
\end{aligned}
\end{equation}
where $g^p(\boldsymbol{\theta}) $ is a probability bound for the random variable $g(\boldsymbol{\theta}) $, and $F_{cdf}$ the cumulative distribution function (CDF). 
In this work, we focus on the mean values and lower-upper bounds of state variables and output values as follows:
\begin{equation}\label{eq12}
\begin{aligned}
\boldsymbol{\mu}_{\boldsymbol {y'}_{t}}&=\mathbb{E}({\boldsymbol {y'}_{t}})\\
\boldsymbol{\alpha}^{lo}_{\boldsymbol {y'}_{t}}&=F^{-1}_{cdf,\boldsymbol {y'}_{t}}( \frac{1}{2}\beta) \\
\boldsymbol{\alpha}^{up}_{\boldsymbol {y'}_{t}}&=F^{-1}_{cdf,\boldsymbol {y'}_{t}}(1- \frac{1}{2}\beta) \\
\boldsymbol{\mu}_{\boldsymbol {z}_{t}}&=\mathbb{E}({\boldsymbol {z}_{t}})\\
\boldsymbol{\alpha}^{lo}_{\boldsymbol {z}_{t}}&=F^{-1}_{cdf,\boldsymbol {z}_{t}}( \frac{1}{2}\beta) \\
\boldsymbol{\alpha}^{up}_{\boldsymbol {z}_{t}}&=F^{-1}_{cdf,\boldsymbol {z}_{t}}(1- \frac{1}{2}\beta) \\
\end{aligned}
\end{equation}
where $\boldsymbol{\mu}_{\boldsymbol {y'}_{t}}$, $\boldsymbol{\alpha}^{lo}_{\boldsymbol {y'}_{t}}$, $\boldsymbol{\alpha}^{up}_{\boldsymbol {y'}_{t}}$ are the mean value, lower and upper bounds of the discrete-time state variable  ${\boldsymbol {y'}_{t}}$ respectively, $\boldsymbol{\mu}_{\boldsymbol {z}_{t}}$,$\boldsymbol{\alpha}^{lo}_{\boldsymbol {z}_{t}}$,$\boldsymbol{\alpha}^{up}_{\boldsymbol {z}_{t}}$ are then mean value, lower and upper bounds of discrete-time output values ${\boldsymbol {z}_{t}}$ respectively, $F_{\boldsymbol {y'}_{t}}$ and $F_{\boldsymbol {z}_{t}}$ are the CDFs of ${\boldsymbol {y'}_{t}}$ and ${\boldsymbol {z}_{t}}$, $\beta$ is the small parameter value to compute lower and upper bounds 
with a high confidence level. In this work, $\beta$ is set to be 0.05  with a 95\% level of confidence.

In this manner, the above large-scale stochastic nonlinear programming Eq.(\ref{eq3}) can be transformed into the following Eq.(\ref{eq13}):
\begin{equation}\label{eq13}
 \begin{aligned}
\min_{\boldsymbol{\lambda'}_{t}} \quad & G''(\boldsymbol{\mu}_{\boldsymbol {y'}_{t}},\boldsymbol{\alpha}^{lo}_{\boldsymbol {y'}_{t}},\boldsymbol{\alpha}^{up}_{\boldsymbol {y'}_{t}},\boldsymbol{\mu}_{\boldsymbol{z}_{t}},\boldsymbol{\alpha}^{lo}_{\boldsymbol{z}_{t}},\boldsymbol{\alpha}^{up}_{\boldsymbol{z}_{t}},\boldsymbol{\lambda'}_{t})  \\   
s.t.  \ \boldsymbol{\mu}_{\boldsymbol {y'}_{t+1}}&=f_{\boldsymbol{\mu}}(\boldsymbol{\mu}_{\boldsymbol {y'}_{t}},\boldsymbol{\alpha}^{lo}_{\boldsymbol {y'}_{t}},\boldsymbol{\alpha}^{up}_{\boldsymbol{y'}_{t}},\boldsymbol{\lambda'}_{t})+ \boldsymbol{w}_{t}, t=0,1,2....k\\
\ \boldsymbol{\alpha}^{lo}_{\boldsymbol {y'}_{t+1}}&=f_{\boldsymbol{\alpha}^{lo}}(\boldsymbol{\mu}_{\boldsymbol {y'}_{t}},\boldsymbol{\alpha}^{lo}_{\boldsymbol {y'}_{t}},\boldsymbol{\alpha}^{up}_{\boldsymbol{y'}_{t}},\boldsymbol{\lambda'}_{t})+ \boldsymbol{w}_{t}, t=0,1,2....k\\
\ \boldsymbol{\alpha}^{up}_{\boldsymbol {y'}_{t+1}}&=f_{\boldsymbol{\alpha}^{up}}(\boldsymbol{\mu}_{\boldsymbol {y'}_{t}},\boldsymbol{\alpha}^{lo}_{\boldsymbol {y'}_{t}},\boldsymbol{\alpha}^{up}_{\boldsymbol{y'}_{t}},\boldsymbol{\lambda'}_{t})+ \boldsymbol{w}_{t}, t=0,1,2....k\\
\ \boldsymbol{\mu}_{\boldsymbol {z}_{t+1}}&=h_{\boldsymbol{\mu}}(\boldsymbol{\mu}_{\boldsymbol {y'}_{t+1}},\boldsymbol {\alpha}^{lo}_{\boldsymbol {y'}_{t+1}},\boldsymbol {\alpha}^{up}_{\boldsymbol{y'}_{t+1}},\boldsymbol{\lambda'}_{t})+ \boldsymbol{v}_{t}, t=0,1,2....k\\
\ \boldsymbol {\alpha}^{lo}_{\boldsymbol {z}_{t+1}}&=h_{\boldsymbol {\alpha}^{lo}}(\boldsymbol {\mu}_{\boldsymbol {y'}_{t+1}},\boldsymbol {\alpha}^{lo}_{\boldsymbol {y'}_{t+1}},\boldsymbol {\alpha}^{up}_{\boldsymbol{y'}_{t+1}},\boldsymbol{\lambda'}_{t})+\boldsymbol{v}_{t}, t=0,1,2....k\\
\ \boldsymbol {\alpha}^{up}_{\boldsymbol {z}_{t+1}}&=h_{\boldsymbol {\alpha}^{up}}(\boldsymbol {\mu}_{\boldsymbol {y'}_{t+1}},\boldsymbol {\alpha}^{lo}_{\boldsymbol {y'}_{t+1}},\boldsymbol {\alpha}^{up}_{\boldsymbol{y'}_{t+1}},\boldsymbol{\lambda'}_{t})+ \boldsymbol{v}_{t}, t=0,1,2....k\\
\ \boldsymbol {\mu}_{\boldsymbol {y'}_{0}}&=\boldsymbol{y}_{0} \\
 g''_{cons}&(\boldsymbol {\mu}_{\boldsymbol {y'}_{t}},\boldsymbol {\alpha}^{lo}_{\boldsymbol {y'}_{t}},\boldsymbol {\alpha}^{up}_{\boldsymbol {y'}_{t}},\boldsymbol {\mu}_{\boldsymbol{z}_{t}},\boldsymbol {\alpha}^{lo}_{\boldsymbol{z}_{t}},\boldsymbol {\alpha}^{up}_{\boldsymbol{z}_{t}},\boldsymbol{\lambda'}_{t}) \leq 0, t=0,1,2....k+1\\
\end{aligned}
\end{equation}where $G''$ the deterministic objective formulation of $G'$ through statistical moments and bounds $\boldsymbol{\mu}_{\boldsymbol {y'}_{t}}, \boldsymbol{\alpha}^{lo}_{\boldsymbol {y'}_{t}}, \boldsymbol{\alpha}^{up}_{\boldsymbol {y'}_{t}}, \boldsymbol{\mu}_{\boldsymbol{z}_{t}}, \boldsymbol{\alpha}^{lo}_{\boldsymbol{z}_{t}}, \boldsymbol{\alpha}^{up}_{\boldsymbol{z}_{t}}$, $f_{\mu}, f_{\alpha^{lo}}, f_{\alpha^{up}}, h_{\mu}, h_{\alpha^{lo}}, h_{\alpha^{up}}$ and $g''_{cons}$ are transformed formulations in similar manners.

\subsection{Proper orthogonal decomposition}
Although the fast PCE-based uncertainty propagation method decreases the computational complexity of calculating statistical moments and bounds, the reduced problem as illustrated in Eq.(\ref{eq13}), is still black-box, high-dimensional and non-convex. %
Double model techniques, POD and RNNs, would be employed to capture the high-dimensional dynamics of statistical moments and bounds (state variables and output values $\boldsymbol{\mu}_{\boldsymbol {y'}_{t}}, \boldsymbol{\alpha}^{lo}_{\boldsymbol {y'}_{t}}, \boldsymbol{\alpha}^{up}_{\boldsymbol{y'}_{t}}, \boldsymbol{\mu}_{\boldsymbol{z}_{t}}, \boldsymbol{\alpha}^{lo}_{\boldsymbol{z}_{t}}, \boldsymbol{\alpha}^{up}_{\boldsymbol{z}_{t}}$) obtained from the fast PCE method. The key idea of the POD method is to project the high-dimensional dynamics of the computed statistical moments/bounds onto a low-dimensional subspace, which requires efficient sampling methods in order to collect enough representative snapshots for a range of parameter values. Details about the theory and applications of POD can be found in the literature \cite{liang2002proper}. 

In this work, Latin hypercube (LHC) sampling on the space of design variables $\boldsymbol{\lambda'}$, combining PCE methods for discretised nonlinear dynamic system Eq.(\ref{eq2}) , was used to get high-dimensional dynamic data sets of statistics$(\boldsymbol {D} \in \mathbb{R}^{N_{\boldsymbol{\lambda'}} \times N} ,[Y_{\boldsymbol{\mu}_{\boldsymbol {y} }} \in \mathbb{R}^{m\times N'_{y'}},Y_{\boldsymbol{\alpha}^{lo}_{\boldsymbol {y} }}\in \mathbb{R}^{m\times N'_{y'} },Y_{\boldsymbol{\alpha}^{up}_{\boldsymbol {y} }}\in \mathbb{R}^{m \times N'_{y'}},Y_{\boldsymbol{\mu}_{\boldsymbol {z} }}\in \mathbb{R}^{m \times N'_{z}},Y_{\boldsymbol{\alpha}^{lo}_{\boldsymbol {z} }}\in \mathbb{R}^{m \times N'_{z}},Y_{\boldsymbol{\alpha}^{up}_{\boldsymbol {z} }}\in \mathbb{R}^{m \times N'_{z}}])$. Here $m\in \mathbb{N}$ is the number of discrete interval points, a large number for distributed parameter systems,  $N\in \mathbb{N}$ is the number of LHC samples. $N'_{y'}=N_{y'}*N*k$ is the number of discrete time-space points for discretised variable $\boldsymbol{y'}$ while $N'_{z}=N_{z}*N*k$ denotes the number of discrete time-space points for discretised variable $\boldsymbol{z}$,  $Y_{\boldsymbol{\mu}_{\boldsymbol {y} }}, Y_{\boldsymbol{\alpha}^{lo}_{\boldsymbol {y} }}, Y_{\boldsymbol{\alpha}^{up}_{\boldsymbol {y} }}, Y_{\boldsymbol{\mu}_{\boldsymbol {z} }}, Y_{\boldsymbol{\alpha}^{lo}_{\boldsymbol {z} }}
, Y_{\boldsymbol{\alpha}^{up}_{\boldsymbol {z} }}$ are the high-dimensional dynamic data sets of statistical moments and bounds ${\boldsymbol{\mu}_{\boldsymbol {y} }}, {\boldsymbol{\alpha}^{lo}_{\boldsymbol {y} }},{\boldsymbol{\alpha}^{up}_{\boldsymbol {y} }}, {\boldsymbol{\mu}_{\boldsymbol {z} }},   {\boldsymbol{\alpha}^{lo}_{\boldsymbol {z} }}
, {\boldsymbol{\alpha}^{up}_{\boldsymbol {z} }}$ respectively through the methods of PCE and snapshots. 
LHC method can systematically generate samples, covering the whole design space and maximizing the difference among the produced samples. Given a relatively large number of samples, LHC strategy can  fill the design space to represent the dynamic features of complex systems. More details about sampling techniques and discussions can be found in the literature \cite{tao2021model}.  

For the high-dimensional dynamic data set $Y_{\mu_{\boldsymbol {y} }}$ over a finite spatial interval $\Omega' \in \mathbb{R}$, the POD method aims to choose a "small" set of dominant POD dynamic modes 
%
$P'_{{\boldsymbol{\mu}_{\boldsymbol {y} }}}= (p_1,p_2,...,p_{a_{{{\boldsymbol{\mu}}_{\boldsymbol {y} }}}})\in \mathbb{R}^{a_{{\mu_{\boldsymbol {y} }}}\times N'_{y'}  }$ ($a_{{\boldsymbol{\mu}_{\boldsymbol {y} }}}\in \mathbb{N}$ is the number of POD modes) through projecting the high-dimensional dynamics of $Y_{\boldsymbol{\mu}_{\boldsymbol {y} }}$ onto the subspace $\mathbb{P'}$ of the $a_{{\boldsymbol{\mu}_{\boldsymbol {y} }}}$POD modes.
\begin{equation}\label{eq14}
\boldsymbol{U_1}=\boldsymbol{P'_{\mu_{\boldsymbol {y} }}Y_{\mu_{\boldsymbol {y} }}}
\end{equation}
Here $\boldsymbol{U}\in \mathbb{R}^{a_{{\mu_{\boldsymbol {y} }}}\times N'_{y'}}$ is the projection of the original dynamic data $\boldsymbol{Y_{\mu_{\boldsymbol {y} }}}$ onto $\mathbb {P'}$ and $\boldsymbol{P'_{{\mu_{\boldsymbol {y} }}}}$ is the orthogonal projector.

In the POD method, the matrix $\mathbb{P'}$ could be constructed by the covariance matrix $\boldsymbol{C_y}\in \mathbb{R}^{N_{y'} \times N_{y'} }$ of the output data $\boldsymbol{Y_{\mu_{\boldsymbol {y} }}}$:
\begin{equation}\label{eq15}
\boldsymbol{C_y}=\frac{1}{N'_{y'} -1}\boldsymbol{Y_{\mu_{\boldsymbol {y} }}Y_{\mu_{\boldsymbol {y} }}^T}
\end{equation}
Here we seek to maximise variance and minimise covariance between data,  i.e. maximising its diagonal elements, while minimising the off-diagonal elements of $\boldsymbol{C_y}$. 

This is equivalent to implementing singular value decomposition (SVD) on $\boldsymbol{C_y}$:  
\begin{equation}\label{eq16}
\boldsymbol{C_y}=\boldsymbol{Z^TZ}=(\frac{1}{\sqrt{N'_{y'}-1}}\boldsymbol{Y_{\mu_{\boldsymbol {y} }}}^T)^T(\frac{1}{\sqrt{N'_{y'}-1}}\boldsymbol{Y_{\mu_{\boldsymbol {y} }}}^T)=\boldsymbol{V_1D_1V_1^T}
\end{equation}
where $D_1 \in \mathbb{R}^{N'_{y'}\times N'_{y'}}$ is a diagonal matrix whose diagonal elements are the eigenvalues of $Z^TZ$ and $V_1$ is the orthogonal matrix whose columns are the eigenvectors of $Z^TZ$, which as can be easily shown are equivalent to the POD modes of $\boldsymbol{Y_{\mu_{\boldsymbol {y} }}}$. In fact, we can keep the first $a_{{\boldsymbol {\mu}_{\boldsymbol {y} }}}$POD modes corresponding to the $a_{{\boldsymbol {\mu}_{\boldsymbol {y} }}}$dominant eigenvalues of $\boldsymbol {C_y}$, where usually $a_{{\boldsymbol {\mu}_{\boldsymbol {y} }}}<<N'_{y'}$, hence $\boldsymbol V \in \mathbb{R}^{N'_{y'}\times a_{{\boldsymbol {\mu}_{\boldsymbol {y} }}}}$ and $D_1$ now contains only the $a_{{\boldsymbol {\mu}_{\boldsymbol {y} }}}$ most dominant eigenvalues of the system, $D_1\in \mathbb{R}^{a_{{\boldsymbol {\mu}_{\boldsymbol {y} }}}\times a_{{\boldsymbol {\mu}_{\boldsymbol {y} }}}}$. We can then set $\boldsymbol{P'_{{\mu_{\boldsymbol {y} }}}}=\boldsymbol{V_1}^T$ and perform data reduction through the projection in Eq.(\ref{eq14}).
The original dynamic data sample, $\boldsymbol{Y_{\mu_{\boldsymbol {y} }}}$ can be reconstructed from the projected data:
\begin{equation}\label{eq17}
\boldsymbol{Y_{\mu_{\boldsymbol {y}}}}=\boldsymbol{{P'_{\mu_{\boldsymbol {y} }}}}^T\boldsymbol{U_1}
\end{equation}
and the inverse projection model (reconstruction model) could be obtained as:
\begin{equation}\label{eq18}
\boldsymbol {\mu}_{\boldsymbol {y'}_{t}}=\boldsymbol{{P'_{\mu_{\boldsymbol {y} }}}}^T{\boldsymbol{u}_{\boldsymbol {\mu}_{\boldsymbol {y} }}}
\end{equation}where $\boldsymbol{u}_{\mu_{\boldsymbol {y} }}$ is the reduced low-dimensional variables.

In the similar manner, the inverse projection model (reconstruction models)
for $\boldsymbol {\alpha}^{lo}_{\boldsymbol {y'}_{t}}, \boldsymbol {\alpha}^{up}_{\boldsymbol{y'}_{t}}, \boldsymbol {\mu}_{\boldsymbol{z}_{t}}, \boldsymbol {\alpha}^{lo}_{\boldsymbol{z}_{t}}, \boldsymbol {\alpha}^{up}_{\boldsymbol{z}_{t}}$ could be computed:
\begin{equation}\label{eq19}
\begin{aligned}
\boldsymbol{\alpha}^{lo}_{\boldsymbol {y'}_{t}}&=
\boldsymbol{{P'_{\alpha^{lo}_{\boldsymbol {y} }}}}^T\boldsymbol{u}_{\boldsymbol{\alpha}^{lo}_{\boldsymbol {y}}}\\
\boldsymbol{\alpha}^{up}_{\boldsymbol {y'}_{t}}&= \boldsymbol{{P'_{\alpha^{up}_{\boldsymbol {y} }}}}^T\boldsymbol{u}_{\boldsymbol{\alpha}^{up}_{\boldsymbol {y} }}\\
\boldsymbol{\mu}_{\boldsymbol {z}_{t}}&=\boldsymbol{{P'_{\mu_{\boldsymbol {z} }}}}^T\boldsymbol{u}_{\boldsymbol{\mu}_{\boldsymbol {z} }}\\
\boldsymbol{\alpha}^{lo}_{\boldsymbol {z}_{t}}&= \boldsymbol{{P'_{\alpha^{lo}_{\boldsymbol {z} }}}}^T\boldsymbol{u_{\alpha^{lo}_{\boldsymbol {z} }}}\\
\boldsymbol{\alpha}^{up}_{\boldsymbol {z}_{t}}&= \boldsymbol{{P'_{\alpha^{up}_{\boldsymbol {z} }}}}^T\boldsymbol{u_{\alpha^{up}_{\boldsymbol {z}}}}
\end{aligned} 
\end{equation}
where $\boldsymbol{{P'_{\alpha^{lo}_{\boldsymbol {y} }}}}, \boldsymbol{{P'_{\alpha^{up}_{\boldsymbol {y} }}}}, \boldsymbol{{P'_{\mu_{\boldsymbol {z} }}}}, \boldsymbol{{P'_{\alpha^{lo}_{\boldsymbol {z} }}}}, \boldsymbol{{P'_{\alpha^{lo}_{\boldsymbol {z} }}}}$ are the corresponding orthogonal projectors to $\boldsymbol {\alpha}^{lo}_{\boldsymbol{y'}_{t}}, \boldsymbol{\alpha}^{up}_{\boldsymbol{y'}_{t}}, \boldsymbol{\mu}_{\boldsymbol{z}_{t}}, \boldsymbol{\alpha}^{lo}_{\boldsymbol{z}_{t}}, \boldsymbol{\alpha}^{up}_{\boldsymbol{z}_{t}}$, respectively while $\boldsymbol{u}_{\boldsymbol{\alpha}^{lo}_{\boldsymbol {y}}}, \boldsymbol{u}_{\boldsymbol{\alpha}^{up}_{\boldsymbol {y} }}, \boldsymbol{u}_{\boldsymbol{\mu}_{\boldsymbol {z} }}, \boldsymbol{u_{\alpha^{lo}_{\boldsymbol {z} }}}, \boldsymbol{u_{\alpha^{up}_{\boldsymbol {z}}}} $ are the corresponding reduced variables. 

\subsection{Recurrent neural network}
Previous POD model strategy has been employed to project the high-dimensional nonlinear dynamics of statistical moments and bounds $\boldsymbol{\mu}_{\boldsymbol {y'}_{t}}, \boldsymbol{\alpha}^{lo}_{\boldsymbol {y'}_{t}}, \boldsymbol{\alpha}^{up}_{\boldsymbol{y'}_{t}}, \boldsymbol{\mu}_{\boldsymbol{z}_{t}}, \boldsymbol{\alpha}^{lo}_{\boldsymbol{z}_{t}}, \boldsymbol{\alpha}^{up}_{\boldsymbol{z}_{t}}$  onto the subspaces of  the low-dimensional dominant dynamics of $\boldsymbol{u}_{\mu_{\boldsymbol {y} }}, \boldsymbol{u}_{\boldsymbol{\alpha}^{lo}_{\boldsymbol {y}}}, \boldsymbol{u}_{\boldsymbol{\alpha}^{up}_{\boldsymbol {y} }}, \boldsymbol{u}_{\boldsymbol{\mu}_{\boldsymbol {z} }}, \boldsymbol{u_{\alpha^{lo}_{\boldsymbol {z} }}}, \boldsymbol{u_{\alpha^{up}_{\boldsymbol {z}}}} $ through the methods of snapshots. Then recurrent neural networks (RNNs) are adopted to capture the temporal coefficients on the low-dimensional subspace.  In general, RNNs can be depicted as a folded computational graph, illustrated in Fig.(\ref{f1}), which is unfolded into a series of connected nodes \cite{goodfellow2016deep}. The mathematical formulation is given in the equation:    

\begin{figure}
    \centering
    \includegraphics[width=0.4\textwidth]{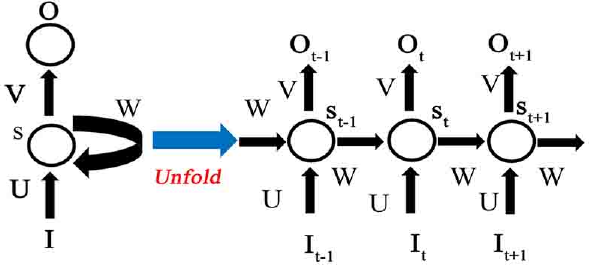}
    \caption{From fold RNN to unfold RNN}
    \label{f1}
\end{figure}
\begin{equation}\label{eq20}
\begin{aligned}
\boldsymbol{s}_{t}&=f_{NN}(\boldsymbol{UI_{t}+Ws_{t-1}+B_1})\\
\boldsymbol{O}_{t}&=g_{NN}(\boldsymbol{Vs_{t}+C_1})
\end{aligned} 
\end{equation}where $s_{t}$ is the time series of the low-dimensional variables $\boldsymbol{u}$ while $\boldsymbol{I}_{t}$ is the time series of the manipulated variables $\boldsymbol{\lambda'}_{t}$, and $\boldsymbol{O}_{t}$ is the output time series. $f_{NN}$ and $g_{NN}$ are RNN functions and output functions with weights $\boldsymbol{U,W,V}$ and biases  $\boldsymbol{B_1}$, $\boldsymbol{C_1}$, respectively. Here we choose the purely linear activation function for the output function $g_{NN}$ since this common choice allows large adjustments for output values.
Then multiple layers are set for RNNs because shallow NNs are easier to be over-fitted, 
which also requires large-scale structures (exponentially larger  number of nodes in one hidden layer) to represent the nonlinear dominant dynamics, 
resulting in intensive computations for the upper-level optimisation and control. Furthermore, the traditional NN activation functions such as the sigmoid and the hyperbolic tangent functions are highly non-convex, leading to multiple local optima. Here the continuous piece-wise affine (PWA) activation function ReLU is adopted to reformulate the trained RNN-based optimisation problem into a MILP problem, which can be applied in an open-loop control mode and be solved using the advanced MILP solver CPLEX \cite{ilog2006cplex}. Here ReLU-RNNs are employed to represent the low-dimensional dynamics of $\boldsymbol{u}_{\mu_{\boldsymbol {y} }}, \boldsymbol{u}_{\boldsymbol{\alpha}^{lo}_{\boldsymbol {y}}}, \boldsymbol{u}_{\boldsymbol{\alpha}^{up}_{\boldsymbol {y} }}, \boldsymbol{u}_{\boldsymbol{\mu}_{\boldsymbol {z} }}, \boldsymbol{u_{\alpha^{lo}_{\boldsymbol {z} }}}, \boldsymbol{u_{\alpha^{up}_{\boldsymbol {z}}}} $. 
Deep rectifier NN based MILP problems have been formulated in previous studies \cite{grimstad2019relu}.
However, the combination of uncertainty quantification, model reduction and deep rectifier NNs has not been reported in the literature. The mathematical equations of deep rectifier NNs for $\boldsymbol{u}_{\mu_{\boldsymbol {y} }}=( {u}_{\mu_{\boldsymbol {y} }1},{u}_{\mu_{\boldsymbol {y} }2},...,{u}_{\mu_{\boldsymbol {y} }a_{{\boldsymbol {\mu}_{\boldsymbol {y} }}}}) \in \mathbb{R}^{a_{{\boldsymbol {\mu}_{\boldsymbol {y} }}}}$ can be reformulated into PWA functions through the big-M method \cite{belotti2010disjunctive} as the equation $F_{\boldsymbol{u}_{\mu_{\boldsymbol {y} }}}(\lambda'_{t},\boldsymbol{u}_{\mu_{\boldsymbol {y} }})=0$.

In a similar manner, the trained RNN models for $\boldsymbol{u}_{\boldsymbol{\alpha}^{lo}_{\boldsymbol {y}}}, \boldsymbol{u}_{\boldsymbol{\alpha}^{up}_{\boldsymbol {y} }}, \boldsymbol{u}_{\boldsymbol{\mu}_{\boldsymbol {z} }}, \boldsymbol{u_{\alpha^{lo}_{\boldsymbol {z} }}}, \boldsymbol{u_{\alpha^{up}_{\boldsymbol {z}}}} $ can be computed as follows: 
\begin{equation}\label{eq22}
\begin{aligned}
F_{\boldsymbol{u}_{\boldsymbol {\alpha}^{lo}_{\boldsymbol {y} }}}(\lambda'_{t},\boldsymbol{u}_{{{\boldsymbol {\alpha}^{lo}_{\boldsymbol {y} }}}})&=0\\
F_{\boldsymbol{u}_{\boldsymbol {\alpha}^{up}_{\boldsymbol {y} }}}(\lambda'_{t},\boldsymbol{u}_{{{\boldsymbol {\alpha}^{up}_{\boldsymbol {y} }}}})&=0\\
F_{\boldsymbol{u}_{\boldsymbol {\mu}_{\boldsymbol {z} }}}(\lambda'_{t},\boldsymbol{u}_{{{\boldsymbol {\mu}_{\boldsymbol {z} }}}})&=0\\
F_{\boldsymbol{u}_{{{\boldsymbol {\alpha}^{lo}_{\boldsymbol {z} }}}}}(\lambda'_{t},\boldsymbol{u}_{{{\boldsymbol {\alpha}^{lo}_{\boldsymbol {z}} }}})&=0\\
F_{\boldsymbol{u}_{{{\boldsymbol {\alpha}^{up}_{\boldsymbol {z} }}}}}(\lambda'_{t},\boldsymbol{u}_{{{\boldsymbol {\alpha}^{up}_{\boldsymbol {z}} }}})&=0
\end{aligned} 
\end{equation}
therefore Eq.(\ref{eq13}) can be transformed into the following large-scale deterministic programming problem:
%
\begin{equation}\label{eq23}
 \begin{aligned}
\min_{\boldsymbol{\lambda'}_{t}} \quad & G''(\boldsymbol{\mu}_{\boldsymbol {y'}_{t}},\boldsymbol{\alpha}^{lo}_{\boldsymbol {y'}_{t}},\boldsymbol {\alpha}^{up}_{\boldsymbol {y'}_{t}},\boldsymbol {\mu}_{\boldsymbol{z}_{t}},\boldsymbol {\alpha}^{lo}_{\boldsymbol{z}_{t}},\boldsymbol {\alpha}^{up}_{\boldsymbol{z}_{t}},\boldsymbol{\lambda'}_{t})  \\   
s.t. \ &F_{\boldsymbol{u}_{\mu_{\boldsymbol {y} }}}(\lambda'_{t},\boldsymbol{u}_{\mu_{\boldsymbol {y} }})=0\\
&F_{\boldsymbol{u}_{\boldsymbol {\alpha}^{lo}_{\boldsymbol {y} }}}(\lambda'_{t},\boldsymbol{u}_{{{\boldsymbol {\alpha}^{lo}_{\boldsymbol {y} }}}})=0\\
&F_{\boldsymbol{u}_{\boldsymbol {\alpha}^{up}_{\boldsymbol {y} }}}(\lambda'_{t},\boldsymbol{u}_{{{\boldsymbol {\alpha}^{up}_{\boldsymbol {y} }}}})=0\\
&F_{\boldsymbol{u}_{\boldsymbol {\mu}_{\boldsymbol {z} }}}(\lambda'_{t},\boldsymbol{u}_{{{\boldsymbol {\mu}_{\boldsymbol {z} }}}})=0\\
&F_{\boldsymbol{u}_{{{\boldsymbol {\alpha}^{lo}_{\boldsymbol {z} }}}}}(\lambda'_{t},\boldsymbol{u}_{{{\boldsymbol {\alpha}^{lo}_{\boldsymbol {z}} }}})=0\\
&F_{\boldsymbol{u}_{{{\boldsymbol {\alpha}^{up}_{\boldsymbol {z} }}}}}(\lambda'_{t},\boldsymbol{u}_{{{\boldsymbol {\alpha}^{up}_{\boldsymbol {z}} }}})=0\\
 g''_{cons}&(
 \boldsymbol{{P'_{\mu_{\boldsymbol {y} }}}}^T{\boldsymbol{u}_{\boldsymbol {\mu}_{\boldsymbol {y} }}},
 \boldsymbol{{P'_{\alpha^{lo}_{\boldsymbol {y} }}}}^T\boldsymbol{u}_{\boldsymbol{\alpha}^{lo}_{\boldsymbol {y}}},
 \boldsymbol{{P'_{\alpha^{up}_{\boldsymbol {y} }}}}^T\boldsymbol{u}_{\boldsymbol{\alpha}^{up}_{\boldsymbol {y} }},\\
 &\boldsymbol{{P'_{\mu_{\boldsymbol {z} }}}}^T\boldsymbol{u}_{\boldsymbol{\mu}_{\boldsymbol {z} }},
 \boldsymbol{{P'_{\alpha^{up}_{\boldsymbol {z} }}}}^T\boldsymbol{u_{\alpha^{up}_{\boldsymbol {z}}}},
 \boldsymbol{{P'_{\alpha^{lo}_{\boldsymbol {z} }}}}^T\boldsymbol{u_{\alpha^{lo}_{\boldsymbol {z} }}},
 \boldsymbol{\lambda'}_{t}) \leq 0, \\
 t&=0,1,2....k+1
\end{aligned}
\end{equation}so far, reduced surrogate models are built to represent the high-dimensional dynamics of statistical moments and/or bounds for large-scale distributed parameter systems. The detailed PCE-POD-RNN model construction algorithm is summarised in Algorithm 1 below:
 \begin{algorithm}
 \caption{PCE-POD-RNN model construction algorithm}
 \begin{algorithmic}[1]
 \renewcommand{\algorithmicrequire}{\textbf{Inputs:}}
 \renewcommand{\algorithmicensure}{\textbf{Outputs:}}
 \REQUIRE $N_1$, $N_2$, $\boldsymbol{w}_{t}$, $\boldsymbol{v}_{t}$
\ENSURE  $\boldsymbol{P'_{{\mu_{\boldsymbol {y} }}}}, \boldsymbol{{P'_{\alpha^{lo}_{\boldsymbol {y} }}}}, \boldsymbol{{P'_{\alpha^{up}_{\boldsymbol {y} }}}}, \boldsymbol{{P'_{\mu_{\boldsymbol {z} }}}}, \boldsymbol{{P'_{\alpha^{lo}_{\boldsymbol {z} }}}}, \boldsymbol{{P'_{\alpha^{up}_{\boldsymbol {z} }}}}$, $F_{\boldsymbol{u}_{\mu_{\boldsymbol {y} }}}, F_{\boldsymbol{u}_{\boldsymbol {\alpha}^{lo}_{\boldsymbol {y} }}}, F_{\boldsymbol{u}_{\boldsymbol {\alpha}^{up}_{\boldsymbol {y} }}}, F_{\boldsymbol{u}_{\boldsymbol {\mu}_{\boldsymbol {z} }}}, F_{\boldsymbol{u}_{{{\boldsymbol {\alpha}^{lo}_{\boldsymbol {z} }}}}},  F_{\boldsymbol{u}_{{{\boldsymbol {\alpha}^{up}_{\boldsymbol {z} }}}}}$
\STATE Generate $N_1$ Latin hypercube samples of $\lambda$ and collect the distributed dynamic trajectories with system and output noises $\boldsymbol{w}_{t}$ and $\boldsymbol{v}_{t}$ 
 \STATE Compute model coefficients of PCE method for time-space trajectories as Eq.(\ref{eq9})
  \STATE Generate $N_2$ Monte Carlo samples (enough) of  uncertain parameters $\boldsymbol{P}$
  \STATE Compute distributed statistical moments and/or bounds (quantities of interest) through the computed PCE models as Eq.(\ref{eq10}-\ref{eq12})
  \STATE Compute each POD global projectors for every statistical moment and/or bound as Eq.(\ref{eq14}-\ref{eq19})
  \STATE Compute each RNN surrogate in the computed low-dimensional POD subspace as Eq.(\ref{eq20})      
    \STATE Check the accuracy of the generated POD-RNN models using fresh samples and dynamic trajectories
	\IF{Model is accurate} 
	\STATE Go to Step 13
	\ELSE 
     \STATE Go back to Step 1 and add more samples or regenerate LHC samples
     \ENDIF
 \RETURN POD projectors and RNN models
 \end{algorithmic} 
 \end{algorithm}
\subsection{Nonlinear model predictive control }
Eq.(\ref{eq23}) gives the open-loop formulation of model reduction based OCP problems, which could be used as an important part of NMPC-based scheme because the solutions of closed-loop NMPC would mainly depend on the online solution of the open-loop problem Eq.({\ref{eq23}}). Here we assume only box constraints for $g^{''}_{cons}$, which commonly include the bounds for manipulated variables and statistics of state and measurement variables as follows:
\begin{equation}\label{eq24}
 \begin{aligned}
\underline{\boldsymbol{\mu}_{\boldsymbol {y'}}}  \leq  \boldsymbol{{P'_{\mu_{\boldsymbol {y} }}}}^T{\boldsymbol{u}_{\boldsymbol {\mu}_{\boldsymbol {y} }}}& \leq \overline{\boldsymbol{\mu}_{\boldsymbol {y'}}}\\ 
\underline{\boldsymbol{\alpha}^{lo}_{\boldsymbol {y'}}} \leq \boldsymbol{{P'_{\alpha^{lo}_{\boldsymbol {y} }}}}^T\boldsymbol{u}_{\boldsymbol{\alpha}^{lo}_{\boldsymbol {y}}}& \leq \overline{\boldsymbol{\alpha}^{lo}_{\boldsymbol {y'}}} \\
\underline{\boldsymbol{\alpha}^{up}_{\boldsymbol {y'}}} \leq\boldsymbol{{P'_{\alpha^{up}_{\boldsymbol {y} }}}}^T\boldsymbol{u}_{\boldsymbol{\alpha}^{up}_{\boldsymbol {y} }}&\leq\overline{\boldsymbol{\alpha}^{up}_{\boldsymbol {y'}}}   \\
\underline{\boldsymbol{\mu}_{\boldsymbol {z}}} \leq   \boldsymbol{{P'_{\mu_{\boldsymbol {z} }}}}^T\boldsymbol{u}_{\boldsymbol{\mu}_{\boldsymbol {z} }}& \leq \overline{\boldsymbol{\mu}_{\boldsymbol {z}}} \\
\underline{\boldsymbol{\alpha}^{lo}_{\boldsymbol {z}}} \leq \boldsymbol{{P'_{\alpha^{lo}_{\boldsymbol {z} }}}}^T\boldsymbol{u_{\alpha^{lo}_{\boldsymbol {z} }}}& \leq \overline{\boldsymbol{\alpha}^{lo}_{\boldsymbol {z}}}   \\
\underline{\boldsymbol{\alpha}^{up}_{\boldsymbol {z}}} \leq \boldsymbol{{P'_{\alpha^{up}_{\boldsymbol {z} }}}}^T\boldsymbol{u_{\alpha^{up}_{\boldsymbol {z}}}}& \leq \overline{\boldsymbol{\alpha}^{up}_{\boldsymbol {z}}}  \\
\underline{\boldsymbol{\lambda'}} \leq {\boldsymbol{\lambda'}}  & \leq \overline{\boldsymbol{\lambda'}}   \\
 &t=0,1,2....k+1
\end{aligned}
\end{equation}
where $\underline{\boldsymbol{\mu}_{\boldsymbol {y'}}}, \underline{\boldsymbol{\alpha}^{up}_{\boldsymbol {y'}}}, \underline{\boldsymbol{\mu}_{\boldsymbol {z}}}, \underline{\boldsymbol{\alpha}^{up}_{\boldsymbol {z}}}, \underline{\boldsymbol{\alpha}^{up}_{\boldsymbol {z}}}, \underline{\boldsymbol{\lambda'}},    \overline{\boldsymbol{\mu}_{\boldsymbol {y'}}}, \overline{\boldsymbol{\alpha}^{lo}_{\boldsymbol {y'}}},  \overline{\boldsymbol{\alpha}^{up}_{\boldsymbol {y'}}}, \overline{\boldsymbol{\mu}_{\boldsymbol {z}}},  \overline{\boldsymbol{\alpha}^{lo}_{\boldsymbol {z}}}, \overline{\boldsymbol{\alpha}^{up}_{\boldsymbol {z}}}, \overline{\boldsymbol{\lambda'}} $ are lower and upper bounds of state variables and manipulated variables, respectively. 

\begin{figure}
    \centering
    \includegraphics[width=0.45\textwidth]{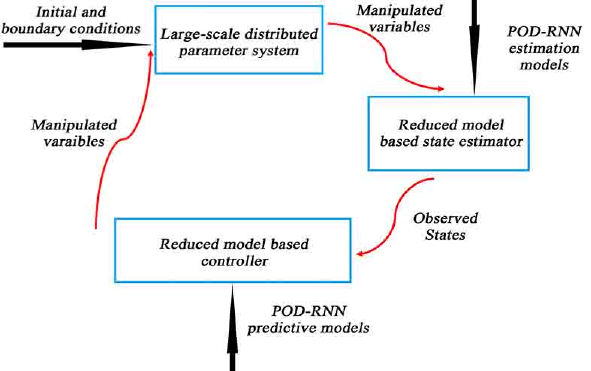}
    \caption{Block diagram of reduced models based observer and NMPC structure }
    \label{f2}
\end{figure}

The objective function $G^{''}$generally aims to minimize the quadratic functions of desired set-points and/or maximize the quantity of interest. In the following case study, we consider maximizing the exit concentration(s) of  products, leading to a shrinking-horizon NMPC.  In this way, the general model reduction-based open loop problems Eq.(\ref{eq23}) can be reformulated into large-scale MILP problems,
solved by the advanced MILP solver CPLEX 12.0. 
Moreover, multiple validated RNN models can be employed to estimate the current statistical moments or bounds of the state variables and quantities of interest, leading to convenient online state estimation using off-line computational models. Then a closed-loop scheme is implemented as shown in  Fig.(\ref{f2}) . 
\section{Case studies}
In this section, we use a tubular chemical reactor and a two-phase  packed bed bioreactor to validate our PCE-POD-RNN-NMPC framework. 
\subsection{Chemical tubular reactor}
This is a chemical engineering application: a tubular reactor as illustrated in Fig.(\ref{f3}), where an exothermic reaction takes place. The reactor model consists of 2 differential equations in dimensionless form as follows: 
\begin{figure}
    \centering
    \includegraphics[width=0.3\textwidth]{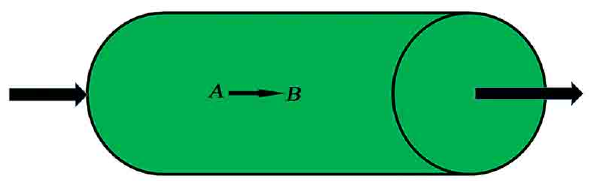}
    \caption{An exothermic tubular reactor with reaction A$\rightarrow$ B}
    \label{f3}
\end{figure}
\begin{equation}\label{eq25}
\begin{aligned}
&\frac{\partial C}{\partial t}=\frac{1}{Pe_1}\frac{\partial^{2}C}{\partial {y1}^{2}}-\frac{\partial C}{\partial {y1}}+D_a(1-C)exp(T/(1+T/\gamma)) \\
&\frac{\partial C}{\partial t}=\frac{1}{LePe_2}\frac{\partial^{2}T}{\partial {y1}^{2}}-\frac{1}{Le}\frac{\partial T}{\partial {y1}}-\frac{\beta}{Le}T+\\
&BD_a(1-C)exp(T/(1+T/\gamma))+\frac{\beta}{Le}T_w   \\  
b.c. &\\
&\frac{\partial C}{\partial {y1}}-{Pe_1}{C} =0, \frac{\partial T}{\partial {y1}}-{Pe_2}{T} =0, \quad at \quad {y1}=0\\
&\frac{\partial C}{\partial {y1}} =0, \frac{\partial T}{\partial {y1}} =0, \quad at \quad {y1}=1
\end{aligned}
\end{equation}Here $C$ and $T$ are the dimensionless concentration and temperature, respectively. $D_a$ denotes the Damköhler number, $Le$ is the Lewis number, $Pe_1$ is the Peclet number for mass transport and $Pe_2$ for heat transport, $\beta$ a dimensionless heat transfer coefficient, $C$ is the dimensionless adiabatic temperature rise, $\gamma $ the dimensionless activation energy and $y1$ the dimensionless longitudinal coordinate. The system parameters are $Pe_1$ = 5, $Pe_2=5$, $Le$ = 1, $Da$ = 0.1, $\beta$ = 1.5, $\gamma $  = 10, $B$ = 12;  
$T_w$ is the adiabatic wall temperature of the cooling problem with the cooling zone. 
A simulator of the model in Eq.(\ref{eq25}) was built with additional system and output noises, and solved  through the $pde$ solver in MATLAB with  200 space discretisation nodes and was subsequently used in inputs/outputs (black-box) mode.  Here inputs include the manipulated variable $T_w$ and two uncertain parameters $D_a \sim N (0.08, 0.008^2)$ and $B\sim N (8, 0.8^2)$  while the outputs are 400 distributed concentration/temperature profiles.  The system and outputs noises satisfy $\boldsymbol{w}_{t}\sim N (\boldsymbol{0}, diag(0.00001*\boldsymbol{1}\in \mathbb{R}^{{400}}))$ and $\boldsymbol{v}_{t}\sim N (\boldsymbol{0}, diag(0.000001*\boldsymbol{1}\in \mathbb{R}^{{400}}))$, respectively.  Meanwhile the reporting time sampling interval was 0.4 for the dynamic simulator.

We aim to improve the performance of this chemical production process through an efficient control strategy under uncertainty based on the constructed black-box simulator. Specifically, 
The objective of the controller design is to maximize $\mathbb{E}(C_{exit})$, the mean value of concentration
at the exit, satisfying the upper bound constraints of temperature $T^{up}(y1,t)$ across the whole reactor. Here Latin hypercube sampling ($N_1=20$) was employed to collect enough representative trajectory samples through the black-box simulator.  Second order polynomial chaos and 4000 realisations ($N_2$) of uncertainty distributions were used to compute the high-dimensional dynamics of statistical moments $\mathbb{E}(C(y1,t))$ and bounds $T^{up}(y1,t)$. Then double model reduction, First POD and then RNNs, was employed to generate the reduced models for online NMPC. Here only 2 dominate POD modes could capture  99.8 \% of the system energy for both $\mathbb{E}(C(y1,t))$ and $T^{up}(y1,t)$, which were then represented by the same size RNNs, 2 hidden layers (15 neurons, 15 neurons) ReLU-based RNNs.  All computations were implemented in MATLAB R2019a on a Desktop (Intel Core(TM) i7-8700 CPU 3.2 292 GHz, 16 GB memory, 64-bit operating system, Windows 10).

\begin{figure}
    \centering
    \includegraphics[width=0.5\textwidth]{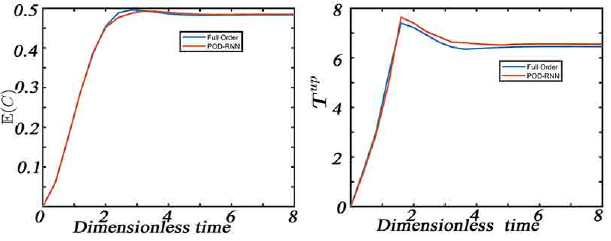}
    \caption{Comparison of full order and POD-RNN time profiles of  $\mathbb{E}(C(y1,t))$ and $T^{up}(y1,t)$ at the exit of the tubular reactor}
    \label{f4}
\end{figure}

\begin{figure}
    \centering
    \includegraphics[width=0.5\textwidth]{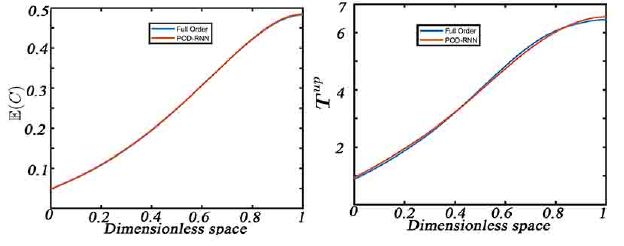}
    \caption{Comparison of full order and POD-RNN space profiles of  $\mathbb{E}(C(y1,t))$ and $T^{up}(y1,t)$ at steady state}
    \label{f5}
\end{figure}

To avoid over-fitting and under-fitting, the data set was randomly divided into a training, a validation and a test set with respective size ratios of 0.70: 0.15: 0.15. The MATLAB Neural Network Toolbox was utilised to fit the weights $\boldsymbol{U,W,V}$ 
and biases  $\boldsymbol{B_1}$, $\boldsymbol{C_1}$
%
by minimizing the mean squared errors (MSE) on the training set using Levenberg-Marquardt algorithm and the early stopping procedure. To obtain the preferable NN structure (numbers of neurons and hidden layers), the training process was repeated until desired accuracy was obtained by adding more nodes and layers. Fig.(\ref{f4}) gives the predicted time profiles of  $\mathbb{E}(C(y1,t))$ and $T^{up}(y1,t)$ at the exit of the tubular reactor while Fig.(\ref{f5}) compares the space profiles at steady state. Both time and space profiles show that the POD-RNN reduced models can approximate the true time-space dynamic process with a high accuracy. 

Then the accurate reduced POD-RNN models could be employed into the model based state estimator and control strategy. Therefore the general open-loop OCP  Eq.(\ref{eq23}) could be applied as formulation Eq.(\ref{eq26}) for the above chemical reaction process, which could be reformulated into MILP problems and iteratively solved by CPLEX 12.0 in GAMS. 
\begin{equation}\label{eq26}
\begin{aligned}
 &\max \limits_{T_{w,t}}\mathbb{E}(C_{exit})\\
s.t.&\ RNN \quad  Equations;\\
&\boldsymbol{{P'_{T^{up}_{\boldsymbol {t} }}}}^T\boldsymbol{u}_{\boldsymbol{T}^{up}_{\boldsymbol {t} }}\leq 4 - \epsilon_1\\
 &  0\leq T_{w,t}\leq 2\\
-\epsilon_2 &\leq  T_{w,t}-T_{w,t-1}\leq \epsilon_2\\
  &  t_n = 1,2...,k \\
\end{aligned}
\end{equation}
Where $\boldsymbol{T}^{up}_{t}(y1)$ is the high-dimensional
discrete upper bounds of distributed temperature at time horizon $t$. $T_{w,t}$ is the discrete manipulated variable, the temperature of cooling zone. $\boldsymbol{u}_{\boldsymbol{T}^{up}_{\boldsymbol {t} }}$ denotes the low-dimensional variables of the projected $\boldsymbol{T}^{up}_{t}(y1) $ through the projector $\boldsymbol{{P'_{T^{up}_{\boldsymbol {t} }}}}$. $\epsilon_1 $ is a relaxation parameter to decrease the impact of approximation errors. $\epsilon_2$ is a limitation parameter for smoothing the jump of manipulated variable. Here $\epsilon_1=\epsilon_2=0.1 $.

\begin{figure}
    \centering
    \includegraphics[width=0.3\textwidth]{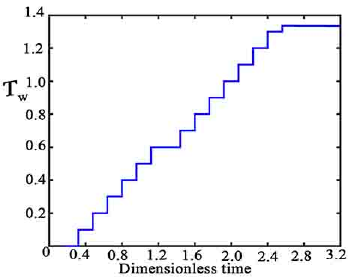}
    \caption{Optimal control policy of heating temperature $T_{w}$}
    \label{f6}
\end{figure}

\begin{figure}
    \centering
    \includegraphics[width=0.5\textwidth]{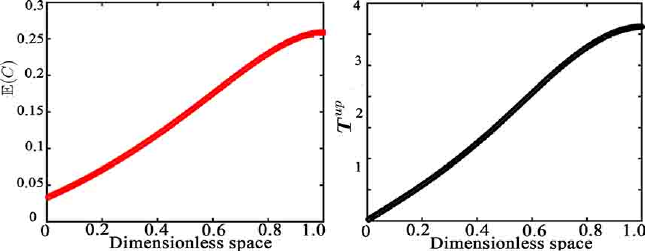}
    \caption{Predictive space profiles of  $\mathbb{E}(C(y1,t))$ and $T^{up}(y1,t)$ at the steady state under the optimal control policy}
    \label{f7}
\end{figure}

\begin{figure}
    \centering
    \includegraphics[width=0.5\textwidth]{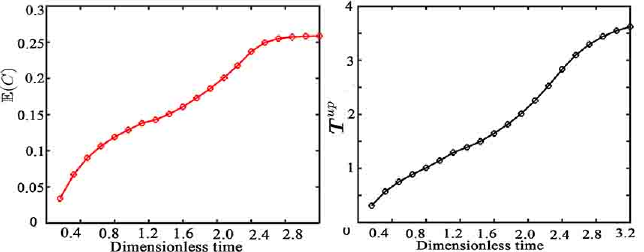}
    \caption{Predictive time profiles of  $\mathbb{E}(C(y1,t))$ and $T^{up}(y1,t)$ at the exit under the optimal control policy}
    \label{f8}
\end{figure}

\begin{figure}
    \centering
    \includegraphics[width=0.5\textwidth]{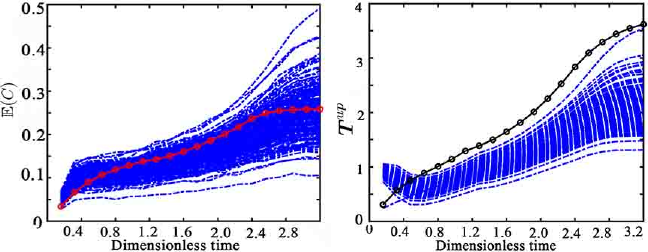}
    \caption{Comparison between random time trajectories (blue dash)and statistical moments (red/black solid) of  $C(y1,t)$ and $T(y1,t)$ at the exit under the optimal control policy}
    \label{f9}
\end{figure}

\begin{figure}
    \centering
    \includegraphics[width=0.5\textwidth]{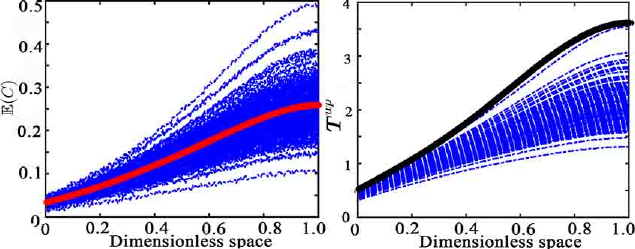}
    \caption{Comparison between random space trajectories (blue dash) and statistical  moments (red/black solid) of  $C(y1,t)$ and $T(y1,t)$ at the steady state under the optimal control policy}
    \label{f10}
\end{figure}

Therefore the reactor system, the reduced model based state estimator and the above open-loop controller make up of the closed-loop control scheme as Fig.(\ref{f2})  for this chemical reaction process. Running the close-loop control system, the close-loop optimal profiles were automatically 
%
generated as the dynamic profiles of the manipulated variable $T_{w}$ in Fig.(\ref{f6}), time dynamic profiles Fig.(\ref{f7}) at the end of the tubular reactor and space profiles Fig.(\ref{f8}) at steady state for both $\mathbb{E}(C(y1,t))$ and $T^{up}(y1,t)$. Then the robustness of the computational framework was tested through running the systems under the computed optimal control policy and 200 random realisations of uncertain parameters. The resulting concentration and temperature profiles are compared in Fig.(\ref{f9}) (time profile) at the exit of the reactor and Fig.(\ref{f10}) (space profile) at steady state. The time profiles of Fig.(\ref{f9}) illustrate that the random time trajectories of  $C(y1,t)$ at the exit  surround the predicted $\mathbb{E}(C(y1,t))$ as expected, while the trajectories of $T(y1,t)$ mostly appear below the predicted $T^{up}(y1,t)$, except for short periods near t=0. One possible reason for this discrepancy is that the errors coming from the model reduction steps are relatively big compared with the small values at the initial stage. 
Meanwhile, the behaviour of spatial profiles is much better.  Across the whole length of reactor, spatially distributed concentrations fluctuate around the predicted $\mathbb{E}(C(y1,t))$ profile. Meanwhile the spatially distributed temperature trajectories always appear below the predicted $T^{up}(y1,t)$. Furthermore, all time-space profiles obey the rigorous constraints $T^{up}(y1,t)=\boldsymbol{{P'_{T^{up}_{\boldsymbol {t} }}}}^T\boldsymbol{u}_{\boldsymbol{T}^{up}_{\boldsymbol {t} }}\leq 4$ for the upper bound variables  $T^{up}(y1,t)$. 

\subsection{Packed-bed bioreactor with immobilised cells}
\begin{figure}
    \centering
    \includegraphics[width=0.15\textwidth]{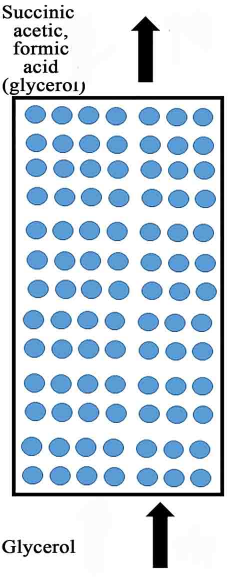}
    \caption{Schematic diagram of packed bed bioreactor with immobilised cells}
    \label{f11}
\end{figure}

Here we extend the model-reduction based robust NMPC strategy to a more complex two-phase (two physical scales) packed bed-bioreactor shown in  Fig.(\ref{f11}). 
Here 
the fermentation process of glycerol to succinic, acetic and formic acid happens in the system of a tubular bioreactor packed with alginate beads, where the biomass is entrapped. The beads with an average diameter of 0.5 cm, were run through multiple fermentations to achieve a constant biomass concentration while the rector was packed until achieving a void coefficient of 0.55 $cm_{packing}^3/cm_{reactor}^3$.  A pump was used to drive the feed through the reactor and the packing. 

Two different control volumes were utilised for the bioprocess: one being the whole bioreactor while the other one being each individaul alginate bead.  Therefore a two-phase model was built for the heterogeneous process, where one-dimensional PDEs Eq.(\ref{eq27}) describe the transport-phenomena in the bulk phase of the reactor, while another set of steady-state PDEs Eq.(\ref{eq28}) illustrate the reaction behaviours in the bead phase \cite{zacharopoulos2021model}.
%

\begin{equation}\label{eq27}
\begin{aligned}
\frac{\partial x_{gly}}{\partial t}&=D_{gly}\frac{\partial^2 x_{gly}}{\partial y1^{2}}-\frac{v}{\epsilon}\frac{\partial x_{gly}}{\partial y1}-R_{tot_{gly}}\\
    \frac{\partial x_{sa}}{\partial t}&=D_{sa}\frac{\partial^2 x_{sa}}{\partial y1^{2}}-\frac{v}{\epsilon}\frac{\partial x_{sa}}{\partial y1}+R_{tot_{sa}}\\
 \frac{\partial x_{aa}}{\partial t}&=D_{aa}\frac{\partial^2 x_{aa}}{\partial y1^{2}}-\frac{v}{\epsilon}\frac{\partial x_{aa}}{\partial y1}+R_{tot_{aa}}\\
 \frac{\partial x_{fa}}{\partial t}&=D_{fa}\frac{\partial^2 x_{fa}}{\partial y1^{2}}-\frac{v}{\epsilon}\frac{\partial x_{fa}}{\partial y1}+R_{tot_{fa}}\\
 x_{gly}\Bigg|_{y1=0}&=x_{gly0}, \frac{\partial x_{gly}}{\partial y}\Bigg|_{y1=L}=0\\
 x_{sa}\Bigg|_{y1=0}&=0, \frac{\partial x_{sa}}{\partial y1}\Bigg|_{y1=L}=0\\
  x_{aa}\Bigg|_{y1=0}&=0, \frac{\partial x_{aa}}{\partial y1}\Bigg|_{y1=L}=0\\
 x_{fa}\Bigg|_{y1=0}&=0, \frac{\partial x_{fa}}{\partial y1}\Bigg|_{y1=L}=0
 \end{aligned}
\end{equation}

$x_{gly}$ ,$x_{sa}$,$x_{aa}$, $x_{fa}$ are the concentrations of glucose, succinic, acetic and formic acid, respectively. $\epsilon$ denotes the void fraction of the packed bed. $v$ is the velocity of the fluid while $R_{tot_{gly}}$, $R_{tot_{sa}}$, $R_{tot_{aa}}$ and $R_{tot_{fa}}$ are the reaction source terms, and $D_{gly}$, $D_{sa}$, $D_{aa}$, $D_{fa}$ are the diffusion coefficients of species in the fluid. \begin{equation}\label{eq28}
\begin{aligned}
    \frac{d}{dr} \Big[ r^2D'_{gly}\frac{\partial x'_{gly}}{\partial r}\Big]-r^2 \cdot r_{gly}&=0\\
  \frac{d}{dr} \Big[ r^2D'_{sa}\frac{\partial x'_{sa}}{\partial r}\Big]+r^2 \cdot r_{sa}=0 \\
\frac{d}{dr} \Big[ r^2D'_{aa}\frac{\partial x'_{aa}}{\partial r}\Big]+r^2 \cdot r_{aa}&=0\\
    \frac{d}{dr} \Big[ r^2D'_{fa}\frac{\partial x'_{fa}}{\partial r}\Big]+r^2 \cdot r_{fa}=0\\
 x'_{gly}\Bigg|_{y1=0}=x_{gly}, \frac{\partial x'_{gly}}{\partial y1}\Bigg|_{y1=L}=0\\
 x'_{sa}\Bigg|_{y1=0}=x_{sa}, \frac{\partial x'_{sa}}{\partial y1}\Bigg|_{y1=L}=0\\
  x'_{aa}\Bigg|_{y1=0}=x_{aa}, \frac{\partial x'_{aa}}{\partial y1}\Bigg|_{y1=L}=0\\
 x'_{fa}\Bigg|_{y1=0}=x_{fa}, \frac{\partial x'_{fa}}{\partial y1}\Bigg|_{y1=L}=0
 \end{aligned}
\end{equation}
Where $D'_{gly}$, $D'_{sa}$, $D'_{aa}$, $D'_{fa}$ are the diffusion coefficients of species in the beads. $r_{gly}$, $r_{sa}$, $r_{aa}$, $r_{fa}$ are the reaction rates of species in the beads.  
\begin{equation}\label{eq29}
\begin{aligned}
R_{tot_{i}}&=\int_{0}^{R} r_{i} \ 4\pi R^2 dR \cdot \rho_{bead}(1-\epsilon),\forall {i} \in \lbrace gly,sa,aa,fa\rbrace \\         
  r_{i}& = (\alpha_{i} \mu + \beta_{i}) \cdot X_{cons}, \forall {i} \in \lbrace gly,sa\rbrace \\
  r_{i}& = \beta_{i} \cdot X_{cons}, \forall {i} \in \lbrace fa,aa\rbrace \\
\mu &=\mu_{max} \cdot \frac{x'_{gly}}{K_{S_{{gly}}}+{x'_{gly}}+\frac{{x'^2_{gly}}}{K_{I_{gly}}}}\cdot \frac{x_{CO_{2}}}{K_{x_{CO_{2}}}+x_{CO_{2}}} \cdot \bigg(1-\frac{x'_{sa}}{{SA}^*}\bigg)^{n_{SA}}
  \end{aligned}
\end{equation}
Where $R$ is the radius of alginate beads, $X_{cons}$ denotes the biomass concentration inside the alginate beads. $\rho_{bead}$ is the density of beads. $\mu$, $\alpha_{i}$ and $\beta_{i}$ denote the coefficients of reaction rates of species. $x_{CO_{2}}$ is the concentration of $CO_{2}$. $K_{I_{gly}}$, $K_{S_{gly}}$,  $K_{x_{CO_{2}}}$,  $SA^{*}$ and $n_{SA}$ are the coefficients of intrinsic kinetics. 

Here $v$=0.1, $L$=20, $D_{gly}$= 0.01, $D_{sa}$=0.01, $D_{aa}$=0.02, $D_{fa}$= 0.02,  $\epsilon$=0.55,  $\rho_{bead}$=2.12,  $R$ =0.15, $X_{cons}$=0.21, $\alpha_{gly}$=2.39, $\alpha_{sa}$=4.5, $\beta_{gly}$=0.187, $\beta_{sa}$=0.21, $\beta_{fa}$=0.011, $\beta_{aa}$=0.0055, $\mu_{max}$=0.2568, $K_{S_{{gly}}}$=5.4,$K_{I_{{gly}}}$=119.99, $x_{CO_{2}}$=0.03, $K_{x_{CO_{2}}}$=0.03, $SA^{*}$ =45.6, $n_{SA}$=5, $D'_{sa}$= 0.00989, $D'_{fa}$= 0.01835, $D'_{aa}$= 0.01384;%

A simulator of the model in Eqs.(\ref{eq27}-\ref{eq29}) was built with additional system and output noises, and solved  through the ODE 113 solver  in MATLAB with 100 space discretisation nodes and was subsequently used in inputs/outputs (black-box) mode. Here inputs include the manipulated variable nominal substrate glycerol concentration  $\overline{x_{gly0}}$ and two uncertain parameters the true substrate glycerol concentration $x_{gly0}\sim U (\overline{x_{gly0}}-
2.5, \overline{x_{gly0}}+2.5)$ and the diffusion coefficients of glycerol in bead phase $D'_{gly}\sim U (0.008, 0.012)$. While the outputs are 400 distributed reactant/product (glycerol, succinic, acetic and formic acid) concentration profiles. The system and outputs noises follow $\boldsymbol{w}_{t}\sim N (\boldsymbol{0}, diag(0.000001*\boldsymbol{1}\in \mathbb{R}^{{400}}))$ and $\boldsymbol{v}_{t}\sim N (\boldsymbol{0}, diag(0.000001*\boldsymbol{1}\in \mathbb{R}^{{400}}))$, respectively. Meanwhile the reporting time sampling interval was 15 hours for the dynamic simulator.

\begin{figure}
    \centering
    \includegraphics[width=0.3\textwidth]{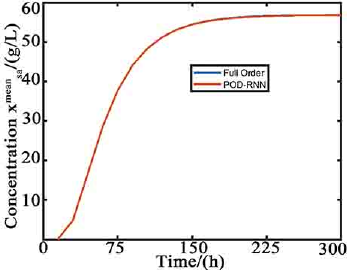}
    \caption{Comparison of expected succinic acid concentration from full systems and predictive POD-RNN models at the exit}
    \label{f12}
\end{figure}

\begin{figure}
    \centering
    \includegraphics[width=0.5\textwidth]{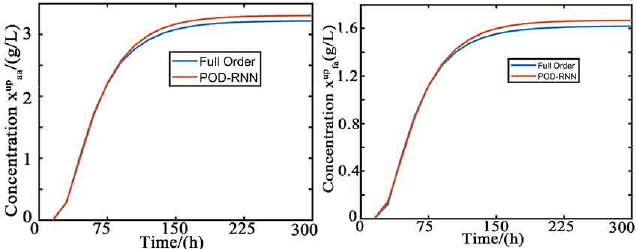}
    \caption{Comparison of upper bound acetic and formic acid concentrations from full systems and predictive POD-RNN models at the exit of the reactor}
    \label{f13}
\end{figure}

\begin{figure}
    \centering
    \includegraphics[width=0.5\textwidth]{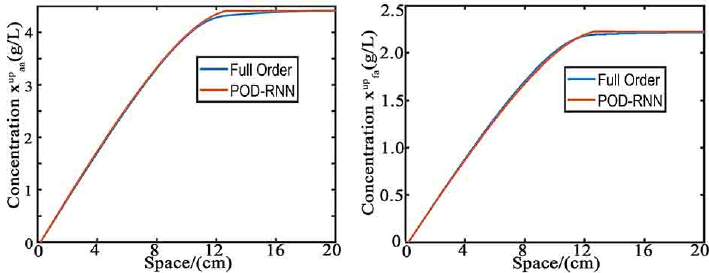}
    \caption{Comparison of upper bound acetic and formic acid concentrations from full systems and predictive POD-RNN models at the steady state}
    \label{f14}
\end{figure}

We aim to improve the performance of this biochemical production process through an efficient control strategy under uncertainty based on the built black-box simulator. 
 Specifically, the objective of controller design is to maximize the expected main product (succinic acid) concentration $\mathbb{E}(x_{sa}^{exit}(t))$) at the exit while observing  upper bound concentration ($x^{up}_{aa}(y1,t)$ and $x^{up}_{fa}(y1,t)$) constraints of  byproducts acetic and formic acid in the whole bioreactor length. Similar to the settings for the previous tubular reactor case study, Latin hypercube sampling ($N_1=20$) to collect enough representative samples using the constructed system simulator in black-box mode. Here second order polynomial chaos and 4000 realisations ($N_2$) of uncertain parameters were used to perform the uncertainty propagation procedure. 
%
%
Then a single RNN (2 hidden layers, 15 neurons and 15 neurons) model was employed to capture the dynamics of the mean value of succinic acid at the exit, 
while POD-RNN double models were built for time -space upper bound profiles of acetic and formic acid. 
Only 2 dominant POD modes are enough to capture 99.8 \% energy of acetic and formic acid dynamics. In addition, two RNNs (2 hidden layers with 15 neurons and 15 neurons, and 2 hidden layers with 20 neurons and 20 neurons) were utilised to construct the low-dimensional dynamic models for acetic and formic acid, respectively. The training process setting was similar to the previous implementation for the chemical tubular reactor.  Then Fig.(\ref{f12}-\ref{f14}) compare the product concentrations from the full system and predictive POD-RNN models. Fig.(\ref{f12}) illustrates that the predicted expected succinic acid concentrations are extremely close to the ones computed from the original system, indicating the high accuracy of the reduced model prediction. In addition, there are only small errors (less than 3 \% on the stable stage) 
between the predictions of the reduced  and the full model as can be seen in Fig.(\ref{f13})-\ref{f14}))h accuracy of the POD-RNN models for the computation of acetic and formic acid profiles.  In general, the constructed reduced surrogate RNN model could accurately predict the complex bioprocess behaviours. 
Then the accurate reduce POD-RNN models could be employed into the model based state estimator and control strategy. Therefore the general open-loop OCP  Eq.(\ref{eq23}) could be applied as formulation Eq.(\ref{eq30}) for the above fermentation process, which could be reformulated into MILP problems.
\begin{equation}\label{eq30}
\begin{aligned}
 &\max \limits_{\overline{x_{gly0}}(t)}\mathbb{E}(x_{sa}^{exit})\\
s.t.&\ RNN \quad \quad  Equations;\\
&\boldsymbol{{P'_{x^{up}_{\boldsymbol {aa} }}}}^T\boldsymbol{u_{x^{up}_{\boldsymbol {aa}}}} \leq 4.1 - \epsilon'_1\\
& \boldsymbol{{P'_{x^{up}_{\boldsymbol {fa} }}}}^T\boldsymbol{u_{x^{up}_{\boldsymbol {fa}}}}  \leq 2.1 - \epsilon'_3\\
 &  50\leq \overline{x_{gly0}}(t)\leq 70\\
-\epsilon'_2 &\leq  \overline{x_{gly0}}(t)-\overline{x_{gly0}}(t-1)\leq \epsilon'_2\\
  &  t = 1,2...,k 
\end{aligned}
\end{equation}
Where  $\epsilon'_1=0.1 $ and  $\epsilon'_3 =0.1$ are the relaxation values to decrease the impact of approximation errors. $\epsilon'_2=1$ is an additional limitation value for smoothing the jump of manipulated variable. $\boldsymbol{u}_{\boldsymbol{x}^{up}_{\boldsymbol {aa} }}$ and  $\boldsymbol{u}_{\boldsymbol{x}^{up}_{\boldsymbol {fa} }}$ are the low dimensional variables of  the projected $\boldsymbol{x}^{up}_{aa} $ and $\boldsymbol{x}^{up}_{fa} $ with the projectors $\boldsymbol{{P'_{x^{up}_{\boldsymbol {aa} }}}}$ and $\boldsymbol{{P'_{x^{up}_{\boldsymbol {fa} }}}}$, respectively.

\begin{figure}
    \centering
    \includegraphics[width=0.5\textwidth]{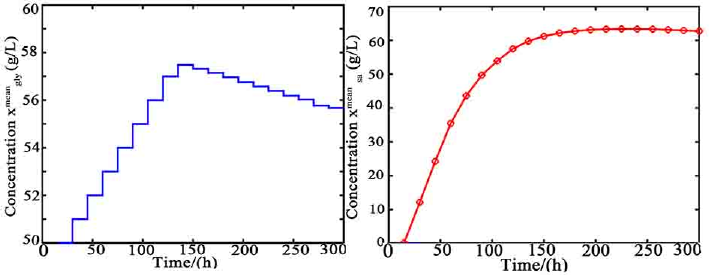}
    \caption{Diagram of the optimal control policy and expected succinic acid concentration  at the exit }
    \label{f15}
\end{figure}

\begin{figure}
    \centering
    \includegraphics[width=0.5\textwidth]{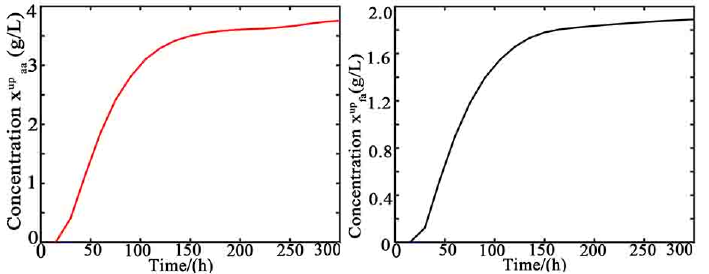}
    \caption{Comparison upper bound acetic and formic acid concentrations at the end of reactor under optimal control policy}
    \label{f16}
\end{figure}

\begin{figure}
    \centering
    \includegraphics[width=0.5\textwidth]{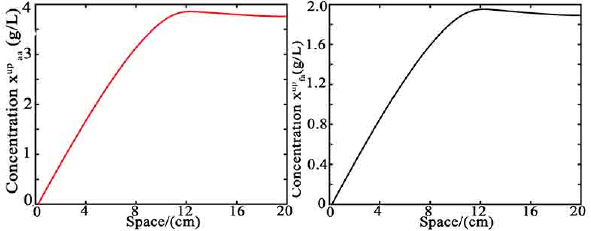}
    \caption{Comparison upper bound acetic and formic acid concentrations at the steady state under optimal control policy}
    \label{f17}
    \par 
\end{figure}

Therefore the reactor system, reduced model based state estimator and above open-loop controller make up of the closed-loop control scheme as Fig.(\ref{f2})  for this biochemical fermentation process. Running the close-loop control system, the close-loop optimal profiles were automatically generated.
Fig.(\ref{f15}) depicts the generated optimal control policy of the manipulated variable substrate glycerol concentration and the dynamic profile of the expected predicted succinic acid at the exit under the control policy. The substrate concentration first went up quickly and then gradually decreased beyond 150 hours, possibly in order to satisfy the requirements for byproduct concentrations. Meanwhile, the computed succinic acid concentration went through initial fast dynamics and then reached steady state. In addition, Fig.(\ref{f16}-\ref{f17}) illustrate the time-space profiles of byproducts acetic and formic acid with similar dynamic behaviours. Time profiles at the exit of the reactor (Fig.(\ref{f16})), showed that the byproduct concentrations grew up quickly, and then became more stable with a little bit final decrease. While space profiles at steady-state stage (Fig.(\ref{f17})) illustrate that the byproduct concentrations increased dramatically along with the reactor and then went down slightly, indicating the nonlinear spatial distributions of  byproduct concentrations. 
 %
 
\begin{figure}
    \centering
    \includegraphics[width=0.3\textwidth]{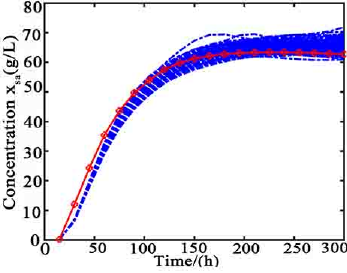}
    \caption{Comparison between random space trajectories (blue dash) and predictive mean value of succinic concentration (red solid)  $x_{sa}$ at the end of reactor under the optimal control policy}
    \label{f18}
\end{figure}

Then the robustness of the proposed control framework was tested through running the systems under 100 random realisations of uncertain parameters within the computed optimal control policy. The resulting succinic acid concentration profile is illustrated in Fig.(\ref{f18}). Generally, the computed concentration profiles are as excepted, with the exception of a small overestimation during the fast dynamic period and a small underestimation for the steady-state value, which may be caused by the model reduction steps. 

\begin{figure}
    \centering
    \includegraphics[width=0.5\textwidth]{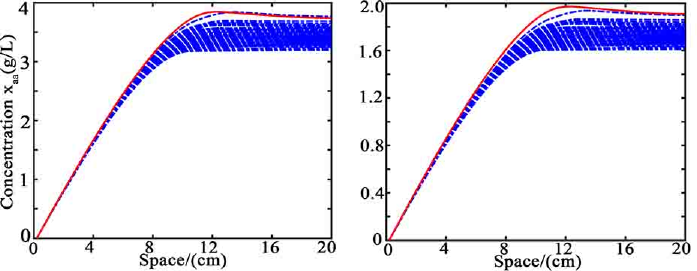}
    \caption{Comparison between random space trajectories (blue dash) and upper bounds (red/black solid) of  $x_{aa}(y,t)$ and $x_{fa}(y,t)$ at the steady state under the optimal control policy}
    \label{f19}
\end{figure}

\begin{figure}
    \centering
    \includegraphics[width=0.5\textwidth]{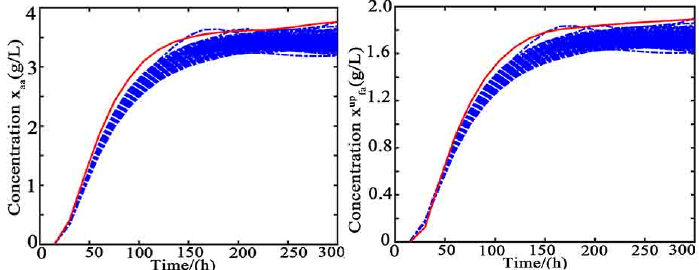}
    \caption{Comparison between random time trajectories (blue dash)and upper bounds (red/black solid) of  $x_{aa}(y,t)$ and $x_{fa}(y,t)$ at the exit under the optimal control policy}
    \label{f20}
\end{figure}

Moreover, the computed acetic and formic acid concentration are displayed in Fig.(\ref{f19}-\ref{f20}).
Almost all of random trajectories followed the dynamic profiles below the predicted upper bounds except for an extremely small initial period for the formic acid and a small period at the end of the exponential phase for both formic and acetic acids. Meanwhile, the gap between the predicted bounds and the random trajectory profiles was not distinct, especially when the profiles were close to the exit of reactor or steady state, which may be caused by the model reductions. Furthermore, all of the 100 random time-space trajectories could satisfy the upper bound constraints $x^{up}_{\boldsymbol {aa} }(y1,t) =\boldsymbol{{P'_{x^{up}_{\boldsymbol {aa} }}}}^T\boldsymbol{u_{x^{up}_{\boldsymbol {aa}}}} \leq 4.1,  x^{up}_{\boldsymbol {fa} }(y1,t)=\boldsymbol{{P'_{x^{up}_{\boldsymbol {fa} }}}}^T\boldsymbol{u_{x^{up}_{\boldsymbol {fa}}}}  \leq 2.1 $, implying the rigorous product quality constraint could be guaranteed with the relation parameters under the computed optimal control policy. 

%
\section{Conclusion} 

This work has proposed a robust model predictive control strategy for large-scale distributed systems under parametric uncertainty. Firstly, the PCE method is used to efficiently compute the high-dimensional dynamics of statistical moments and probabilistic bounds for time-space state variables and measurement outputs. Then, the double model reductions, POD and RNN, are employed to construct simple but accurate predictive models for the nonlinear dynamics of statistical moments and upper/lower bounds. Next the MILP models are used to reformulate online optimisation scheme and solved globally within the NMPC. The proposed methodology is verified by a typical chemical tubular reactor and a packed bed bioreactor with immobilised cells for production. The two cases shows that the proposed NMPC strategy could efficiently improve process production and also satisfy the requirements of process safety  (time-space temperature constraints) and product quality (time-space byproduct constraints). In the future, less conservative strategy would be exploited to further enhance the process performances.  Moreover, online adaptive control was considered when offline samples are not enough.
\bibliographystyle{IEEEtran}
\bibliography{IEEEabrv,IEEEexample}
\end{document}